%% file: brou_20170113_ArXiv2a.tex
\title{Using Brouwer's fixed point theorem\footnote{To appear in \emph{``A Journey through Discrete Mathematics. A Tribute to Ji\v{r}\'{i} Matou\v{s}ek''}, edited by Martin Loebl, Jaroslav Ne\v{s}et\v{r}il and Robin Thomas, due to be published by Springer.}}

\author{%
{\sc Anders Bj\"orner}
\\
   {\footnotesize Department of Mathematics}\\[-1.5mm]
   {\footnotesize Royal Institute of Technology (KTH)}\\[-1.5mm]
   {\footnotesize 100 44 Stockholm, Sweden}
\and
{\sc Ji\v{r}\'{\i} Matou\v{s}ek} 
\\
   {\footnotesize Department of Applied Mathematics and}\\[-1.5mm]
   {\footnotesize Institute of Theoretical Computer Science}\\[-1.5mm]
   {\footnotesize Charles University, Prague, and}\\[-1.5mm]
   {\footnotesize Institute of Theoretical Computer Science}\\[-1.5mm]
   {\footnotesize ETH Zurich, 8092 Zurich, Switzerland}
\and 
{\sc G\"unter M. Ziegler}
\\
   {\footnotesize Institute of Mathematics}\\[-1.5mm]
   {\footnotesize Freie Universit\"at Berlin}\\[-1.5mm]
   {\footnotesize  Arnimallee 2}\\[-1.5mm]
   {\footnotesize  14195 Berlin, Germany}
}

\documentclass[11pt]{article} 
\usepackage[utf8]{inputenc}  
\usepackage[T1]{fontenc}
\usepackage{a4wide}
\usepackage{amssymb,amsmath,amsthm}
\usepackage{epsf,epsfig}
\usepackage{graphicx} 
\usepackage{topmac-brou-arXiv2a}
\usepackage{url,paralist}
\usepackage{color}
\usepackage{epstopdf}
 
\theoremstyle{definition}
\newtheorem{Definition}{Definition}[section]
\theoremstyle{theorem}
\newtheorem{Theorem}[Definition]{Theorem}
\newtheorem{Proposition}[Definition]{Proposition}
\newtheorem{Lemma}[Definition]{Lemma}
\newtheorem{Corollary}[Definition]{Corollary}
\newtheorem{Conjecture}[Definition]{Conjecture}

\theoremstyle{remark}

\newtheorem{Example}[Definition]{Example}
\newtheorem{Examples}[Definition]{Examples}

\date{}
\begin{document}

\maketitle

Brouwer's fixed point theorem from 1911 is a basic result in topology --- 
with a wealth of combinatorial and geometric consequences.
In these lecture notes we present some of them, related to
the game of HEX and to the piercing of multiple intervals.
We also sketch stronger theorems, due to Oliver and others,
and explain their applications to the fascinating (and still not fully 
solved) evasiveness problem.

\begin{small}
	\tableofcontents  
\end{small}

\bigskip 
\setcounter{section}{-1}
\section{Introduction}

The fixed point theorem of Brouwer is one of the most
widely known results of topology.
It says that every continuous map $f:B^d\rightarrow B^d$
of the $d$-dimensional closed unit ball to itself has a fixed point, that is,
a point $x_0\in B^d$ such that $f(x_0)=x_0$.

This result was established by Luitzen Egbertus Jan Brouwer (1881--1960)
at the end of his important 1911 paper \cite{Brouwer11}, in which he also introduced
the fundamental concept (and proof technique) of the mapping degree.
It has many striking and famous applications to problems in
Geometry, Analysis, Game Theory and Combinatorics.

Brouwer's fixed point theorem is in several ways similar to 
the \emph{Borsuk--Ulam theorem} from 1933, which has gotten a lot of 
attention and appreciation for being unusually rich in applications.
For example, the 1978 proofs of the 1955 Kneser conjecture by Lovász and by Bárány
employed the Borsuk--Ulam Theorem  
in order to solve a problem about partitioning a set system,
or equivalently, bounding the chromatic numbers for a certain class of graphs. This
unexpected use of a result from equivariant topology is one of 
the starting points (probably the most famous one)
for the field of ``Topological Combinatorics'' \cite{deLongueville2004}.
We refer to the detailed, elementary exposition in Matou\v{s}ek's book 
``Using the Borsuk--Ulam Theorem'' \cite{Mat-top}.
Current research continues this line of work, using more advanced methods from Equivariant Algebraic Combinatorics;
see for example the text
``Beyond the Borsuk--Ulam Theorem: The Topological Tverberg Story'' \cite{blagojevic2017} in this volume.

In various respects, Brouwer's theorem is a simpler
result than the Borsuk--Ulam theorem: For example, it
is very easy to state (as it does not involve symmetry, or a group action),
and it is quite easy to prove (see below).
It can also easily be derived from the Borsuk--Ulam theorem (see \cite{su1997}),
while indeed it is not as straightforward to obtain ``Borsuk--Ulam from Brouwer.''

Just like the Borsuk--Ulam theorem, Brouwer's theorem 
has many equivalent versions, as well as powerful and useful
extensions. For instance, the Lefschetz fixed point theorem that
works for spaces much more general than a ball,
the Schauder fixed point theorem that works also for compact
balls in infinite-dimensional Banach spaces, the Kakutani fixed point theorem
for set-valued maps, \emph{and so on}. See Shapiro \cite{Shapiro}
for a friendly introduction to fixed point theorems with Analysis applications in mind.

The striking applications of the Brouwer theorem in 
Combinatorics and Geometry seem not to be as well known as the 
applications of the Borsuk--Ulam theorem.
In order to help to remedy this, we present three distinct 
areas of such applications in the three main sections of these lecture notes:
\begin{enumerate}
	\item Brouwer's theorem can be invoked to prove that
		the game of HEX can never end without a winner.
		And indeed, the $d$-dimensional version of this claim
		turns out to be equivalent to Brouwer's theorem!
		This observation of David Gale in his award-winning 1979 paper
		\cite{Gale-hex} may also be counted	
		among the starting points of Topological Combinatorics.
		In our presentation we not only use this to prove the HEX theorem,
		but we also give a combinatorial proof of the HEX theorem
		and derive Brouwer's theorem from this.
	\item Some results about hypergraph matchings and transversals
		have a topological core, to be derived from the Brouwer theorem.
		Our presentation treats one striking instance,
		concerning the relation between packing and transversal numbers 
		for systems of $d$-intervals.
	\item The \emph{Evasiveness conjecture} states that every non-trivial
		monotone graph property is evasive, that is, it does not allow for a
		query strategy that cannot be tricked into checking {\em all} potential edges
		of a graph in order to establish the property.
		This conjecture is still open in general, but the special case of
		a graph on a prime power number of vertices was proved 
		using fixed point theorems of Smith and Oliver. These
		theorems may be seen
		as extensions of Brouwer's. The Appendix to this paper collects
		and sketches the necessary tools.
\end{enumerate}
Further remarkable applications of Brouwer's fixed point theorem
on geometric problems, not treated here, 
include the work by Bondarenko \& Viazovska \cite{BondarenkoViazovska}
on the construction of spherical designs, and the work 
on center points and regression depth by Amenta et al.~\cite{AmentaBernEppsteinTeng}.

Our presentation is based on lecture notes  
that were written about fifteen years ago,
with a history that for some parts goes back nearly thirty years.
These notes can be regarded as a companion or perhaps as a ``prequel'' to Matou\v{s}ek's book \cite{Mat-top}. 

The three main parts do not depend on each other, so they can be
read indepenently.
We refer to \cite{Mat-top} for notation and terminology not explained here.

\bigskip
\section{A game model for Brouwer's fixed point theorem}
\subsection{The Game of HEX}

Let's start with a game: ``HEX'' is a board game for two players,
invented by the ingenious Danish poet, designer and engineer Piet
Hein in 1942 \cite{Gardner}, and rediscovered in 1948 by the mathematician
John Nash \cite{Milnor-nash}, who got a Nobel memorial 
prize in economics in 1994
(for his work on game theory, but not really for this game~\ldots\,).

\def\labelboard{\unitlength 0.9cm
\begin{picture}(0,0)
\put(0.5,3.5){$W$}
\put(3.3,0.2){$B$}
\put(2.8,6.0){$B'$}
\put(5.7,2.6){$W'$}
\end{picture}}

HEX, in Hein's version, 
is played on a rhombical board, as depicted in the figure.
\[
\labelboard
\def\epsfsize#1#2{.36#1}\epsffile{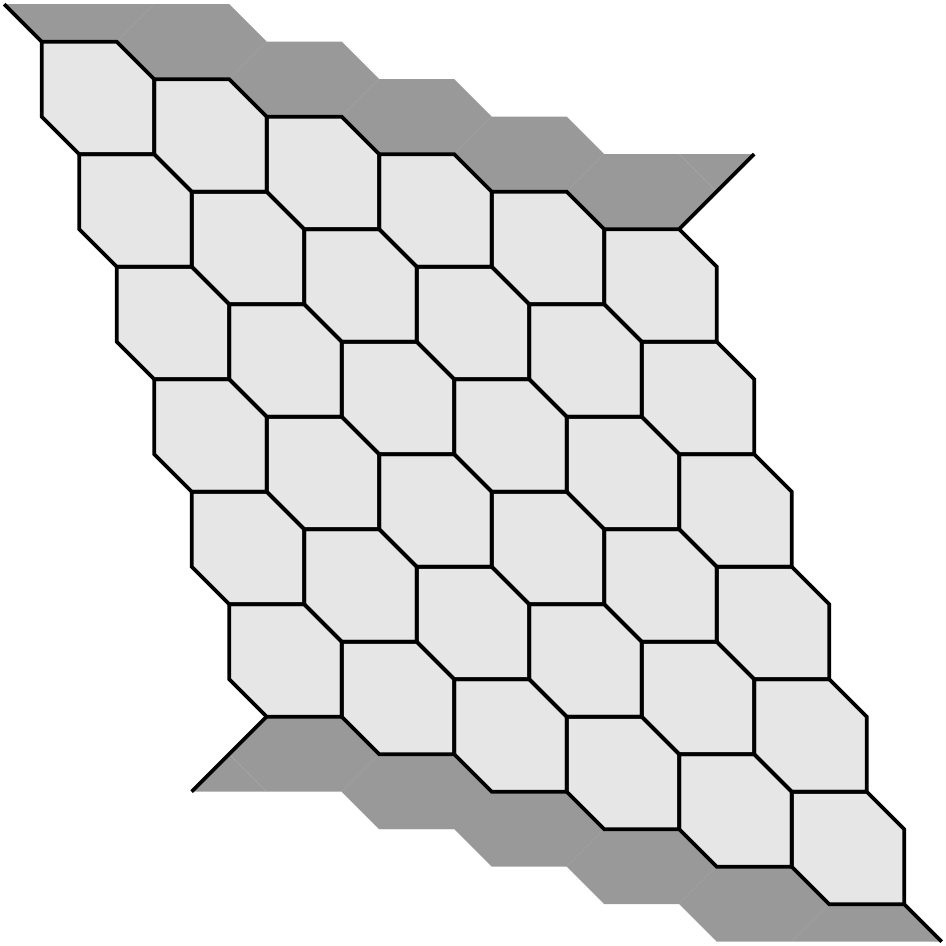}
\]

The rules of the game are simple: There are two players, whom we call
White and Black. The players alternate, with White going first.
Each move consists of coloring one ``grey'' hexagonal tile of the 
board white resp.\ black. 
White has to connect the white borders of
the board (marked $W$ and $W'$) by a path of his white tiles, while 
Black tries to connect $B$ and $B'$ by a black path.
They can't both win: Any winning path for white separates the two
black borders, and conversely. (This isn't hard to prove---however,
the statement is closely related to the Jordan curve theorem, which
is trickier than it may seem when judged at first sight: 
see Exercise~\ref{exer:jct}.)

However, here we concentrate on the opposite statement: There is no
draw possible---when the whole board is covered by black and white
tiles, then there always is a winner. (This is even true if one of
the players has cheated badly and ends up with much more tiles than
his/her opponent! It is also true if the board isn't really ``square,''
that is, if it has sides of unequal lenghts.)
Our next figure depicts a final HEX position---sure enough one of the players
has won, and the proof of the following ``HEX theorem''
will give us a systematic method to find out which one.

\[  
\labelboard
\def\epsfsize#1#2{.36#1}\epsffile{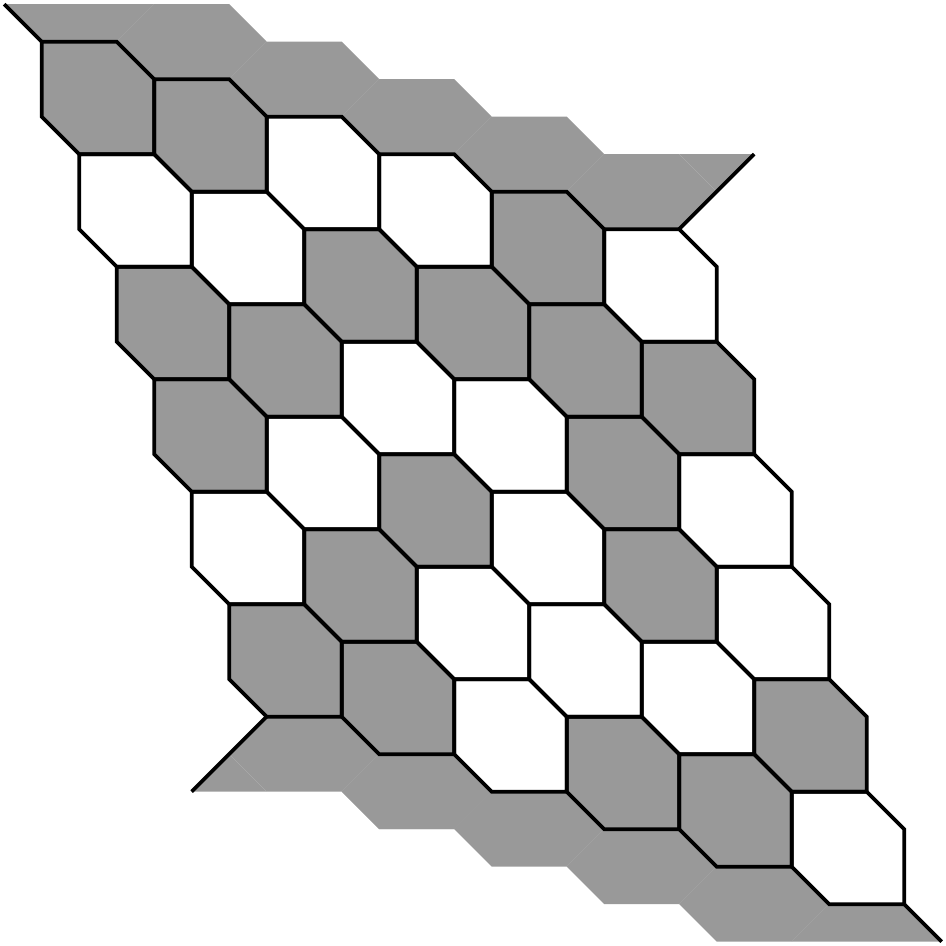}
\]

\begin{Theorem} [The HEX theorem]
If each tile of an $(n \times m)$-HEX board is
colored black or white, then either there is a path of
white tiles that connects the white borders $W$ and $W'$, or there is a
path of black tiles that connects the black borders $B$ and $B'$.
\end{Theorem}

\noindent
Our plan for this section is the following:
\begin{compactitem}[ $\bullet$]
\item We give a simple proof of the HEX theorem.
\item We show that it implies the Brouwer fixed point theorem \ldots
\item \ldots\ and conversely: The Brouwer fixed point theorem implies
the HEX theorem.
\item Then we prove that one of the players has a winning strategy.
\item And then we see that on a square board, the first player can
win, while on an uneven board, the player with the longer borders 
has a strategy to win.
\end{compactitem}
All of this is really quite simple, but it nicely illustrates how a
topological theorem enters the analysis of a discrete situation. 

\proofheader{Proof of the HEX theorem}
For the proof we trace a certain path \emph{between} the black and the white tiles.
It starts
in the lower left-hand corner of the HEX board on the edge that
separates $W$ and $B$. Whenever the path reaches a corner of degree~$3$, 
there will be both colors present at the corner (due to the
edge we reach it from), and so there will be a unique edge to
proceed on that does have different colors on its two sides.

\[
\raisebox{3mm}{\labelboard}
\def\epsfsize#1#2{.36#1}\epsffile{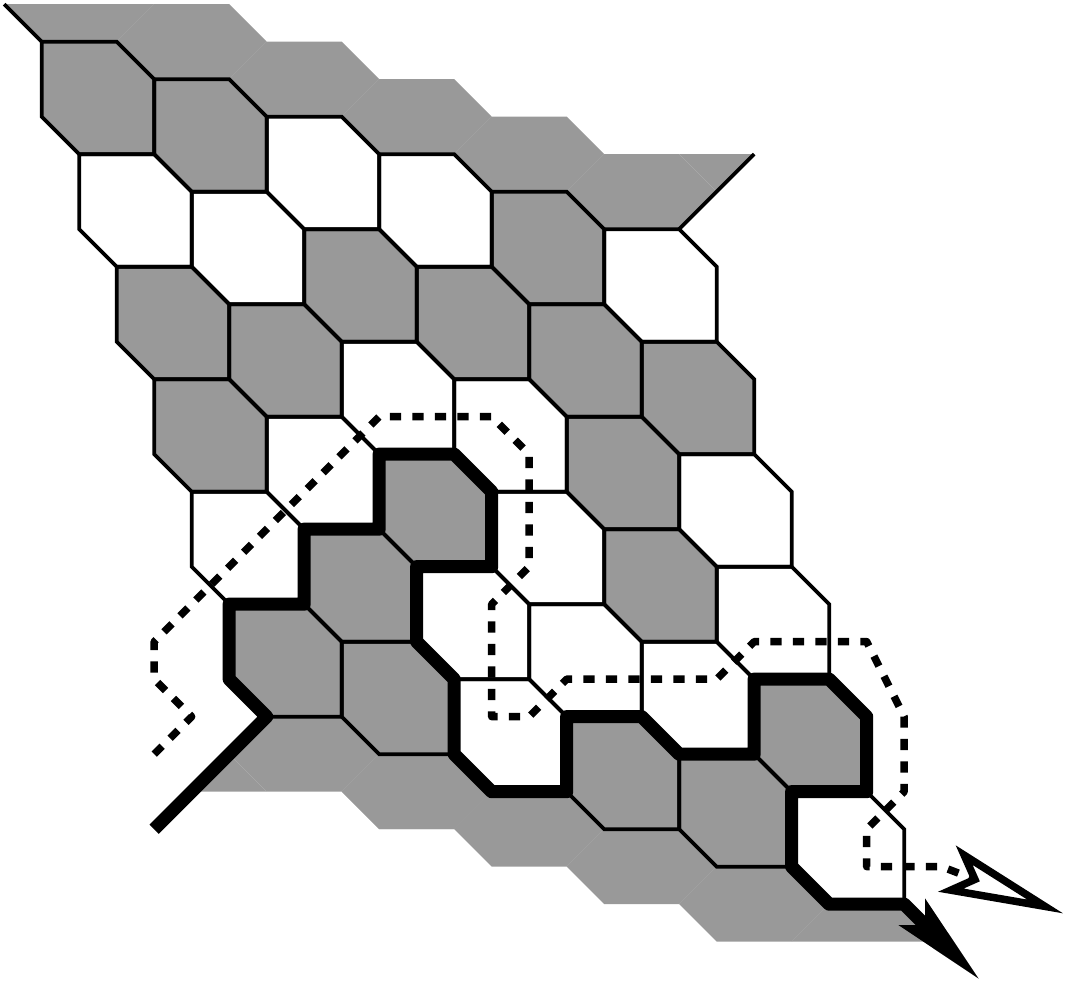}
\]

Our path can never get stuck or branch or turn back onto itself,
otherwise we would have found a vertex
that has one or three edges that separate colors, whereas this
number clearly has to be even at each vertex. Thus the path can be continued
until it leaves the board---that is, until it reaches $W'$ or~$B'$.
But that means that we find a path that connects $W$ to~$W'$, or
$B$ to $B'$, and on its sides keeps a white path of tiles resp.\ a black
path. That is, one of White and Black has won!
\endproof

Now this was easy, and (hopefully) fun. We continue with a
re-interpretation of the HEX board---in Nash's 
version---that buys us two drinks for the price of one: 
\begin{itemize}\itemsep=-2pt
\item[(i)]
a $d$-dimensional version of the HEX theorem, and 
\item[(ii)]
the connection to the Brouwer fixed point theorem.
\end{itemize}

\begin{Definition}[The $d$-dimensional HEX board]
The $d$-dimensional \emph{HEX board} is the graph $H(n,d)$ on the
vertex set $V=\{-1,0,1,\dots,n,n+1\}^d$, in which two vertices
$\vv,\ww\in V$ are connected by an edge if and only if $\vv-\ww\in\{0,1\}^d\cup\{0,-1\}^d$.

The \emph{colors} for the $d$-dimensional HEX game are $1,2,\dots,d$, where we
identify {\rm ``$1 =$\,white''} and {\rm ``$2 =$\,black.''} 
The \emph{interior} of the
HEX board is given by $V'=\{0,1,2,\dots,n\}^d$. All the other vertices,
in $V\setminus V'$, form the \emph{boundary} of the board. The vertices in
the boundary of $H(n,d)$ get {{preassigned colors}} 
\[
\kappa(\vv) = 
\kappa(v_1,\dots,v_d)\assg 
\begin{cases}
    \min\{i\sep v_i=-1 \} &  \textrm{if this exists},\\
    \min\{i\sep v_i=n+1\} &  \textrm{otherwise}.
\end{cases}
\]
\end{Definition}

\[
\def\epsfsize#1#2{.6#1}\epsffile{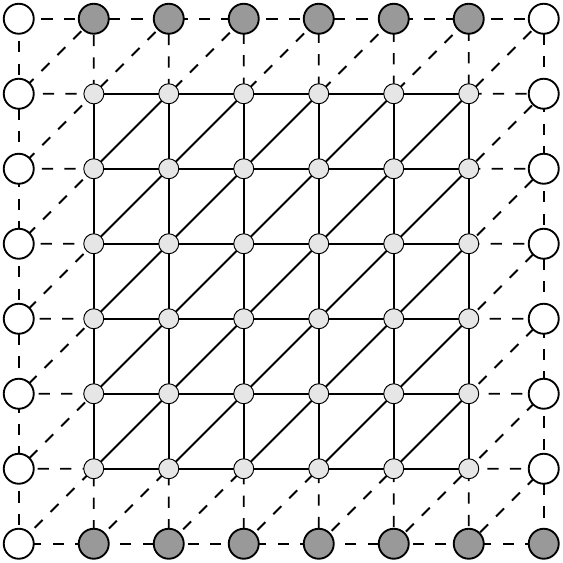}
\]

Our drawing depicts the $2$-dimensional HEX board 
$H(5,2)$, which represents a dual graph for
the $(6\times6)$-board that we
used in our previous figures, with the
preassigned colors on the boundary.

The $d$-dimensional HEX game is  played between $d$ players who
take turns in coloring the interior vertices of~$H(n,d)$. The $i$-th
player \emph{wins} if he\footnote{Using ``he'' here is
not politically correct.} achieves a path of vertices of color $i$ that
connects a vertex whose $i$-th coordinate is $-1$ to a vertex whose
$i$-th coordinate is $n+1$.

\begin{Theorem}[The $d$-dimensional HEX theorem]
For $d$-dimensional HEX at least one of the players reaches his goal: When all interior
vertices of~$H(d,n)$ are colored, then at least one player has won.
\end{Theorem}

\proof
The proof that we used for $2$-dimensional HEX still
works: It just has to be properly translated for the new setting.
For this we first
check that $H(n,d)$ is the graph of a triangulation $\Delta(n,d)$ 
of~$[-1,n+1]^d$, which is given by the \emph{clique complex} 
of~$H(n,d)$. That is, a set of lattice points
$S\subseteq\{-1,0,1,\dots,n+1\}^d$ forms a simplex in
$\Delta(n,d)$ if and only if the points in $S$ are pairwise
connected by edges. To check this, verify that each point
$x\in[{-}1,n{+}1]^d$ lies in the relative interior of a unique simplex,
which is given by
\begin{eqnarray*}
	\lefteqn{\Delta(\xx)\assg 
    \conv\big\{\vv\in\{{-}1,\dots,n+1\}^d\sep}\qquad\qquad\qquad\qquad\\
  & & 
    \lfloor x_i\rfloor\le v_i\le\lceil x_i\rceil
    \mbox{\rm~for all $i$,}\\[1mm]
  &  & 
    \lfloor x_i-x_j\rfloor\le v_i-v_j\le\lceil x_i-x_j\rceil
    \mbox{\rm~for all $i\neq j$}\big\}.
\end{eqnarray*}

Every full-dimensional simplex in $\Delta(n,d)$ has $d+1$
vertices. A simplex $S$ in $\Delta(n,d)$ is
\emph{completely colored} if it has all $d$ colors on its vertices.
Thus each completely colored $d$-simplex in $\Delta$ has exactly
two completely colored facets, which are $(d-1)$-faces of the 
complex $\Delta(n,d)$. Conversely, every completely colored $(d-1)$-face is
contained in exactly two completely colored $d$-simplices---if it is not on the
boundary of~$[-1,n+1]^d$.

With this the (constructive) proof that we gave before for the
$2$-dimensional HEX theorem generalizes to the following: We start
at the $d$-simplex 
\begin{eqnarray*}
	\Delta_0 &\assg&  \conv\{-\one,-\one+\ee_1,-\one+\ee_1+\ee_2,
					\ \dots\ ,-\one+\ee_1+\dots+\ee_{d-1},-\one+\ee_1+\dots+\ee_d\}\\
			 & = &    \conv\{-\one,-\one+\ee_1,-\one+\ee_1+\ee_2,
					\ \dots\ ,-\ee_d,\,\zero\},
\end{eqnarray*}
whose facet ($(d-1)$-face)
$\conv\{-\one,-\one+\ee_1,\dots,-\ee_{d-1}-\ee_d,-\ee_d\}$ is 
completely colored. (Verify this!)
This simplex is shaded in the following figure for $H(5,2)$, which 
depicts the same final position that we considered before.
\[
\def\epsfsize#1#2{.6#1}\epsffile{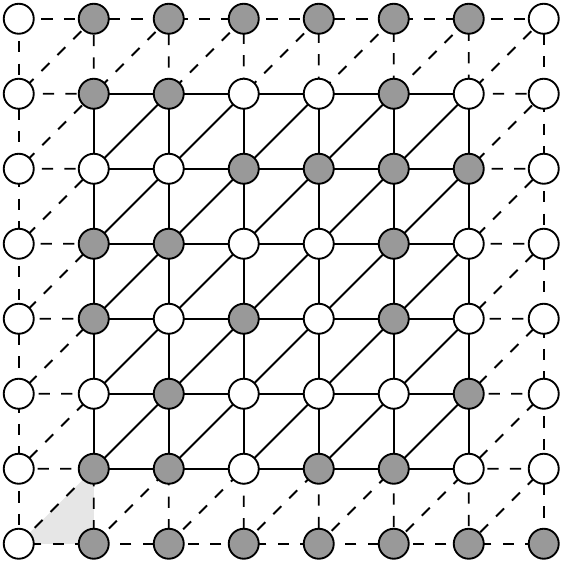}
\]
Now we construct a sequence of completely colored $d$-dimensional
simplices that starts at~$\Delta_0$: We find the 
second completely colored $(d-1)$-face of~$\Delta_0$, find the second
completely colored
$d$-simplex it is contained in, etc. Thus we find a chain of completely colored
$d$-simplices that ends on the boundary of $[-1,n{+}1]^d$---at a
different simplex than the one we started from. In particular, the
last $d$-simplex in the chain has a completely colored facet in the
boundary, and by construction this facet has to lie in a hyperplane
$H^+_i=\{\xx\sep x_i=n+1\}$. 
At this point we check that every completely colored 
$(d-1)$-simplex in the boundary of~$H(n,d)$ is contained in
one of the hyperplanes $H^+_i$, with the sole exception of 
the boundary facet of our starting $d$-simplex.
The chain of $d$-simplices then provides
us with an $i$-colored path from the $i$-colored vertex
\[
-\one + \ee_1+\dots+\ee_{i-1}\in H^-_i=\{\xx\sep x_i=-1\}
\]
to the $i$-colored vertex in $H^+_i$: So the $i$-th player wins.
\endproof

Our drawing illustrates the chain of completely colored simplices
(shaded) and the sequence of (white) vertices for the winning path 
that we get from it.
\[
\def\epsfsize#1#2{.6#1}\epsffile{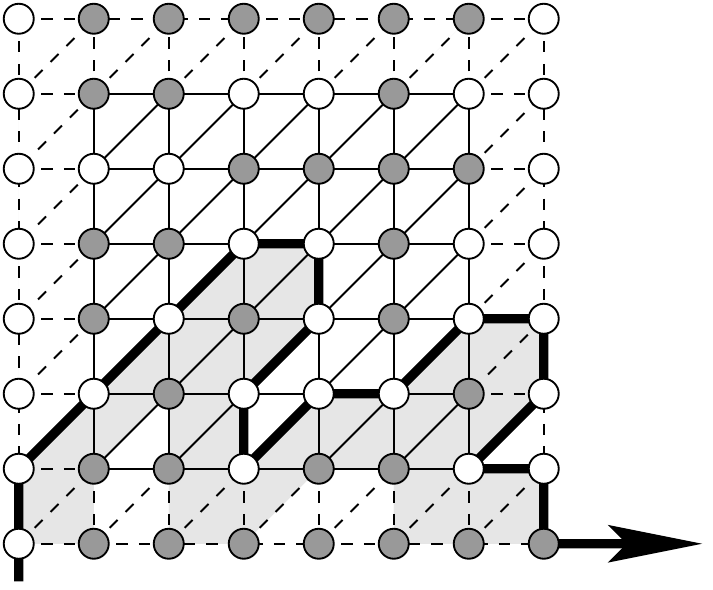}
\]

\vspace{3mm}

\subsection{The Brouwer fixed point theorem}

Now we proceed from the discrete mathematics setting of the 
HEX game to the continuous world of topological fixed point theorems.
Here are three versions of the Brouwer fixed point theorem. 

\begin{Theorem}[Brouwer fixed point theorem]\label{t:brouwer} 
 The following are equivalent (and true):
\begin{compactenum}[ \rm (Br1)]
\item Every continuous map $f\: B^d\to B^d$ has a fixed point.
\item Every continuous map $f\: B^d\to S^{d-1}$ has a fixed point.
\item Every null-homotopic map $f\: S^{d-1}\to S^{d-1}$ has a
fixed point.
\end{compactenum}
\[
\def\epsfsize#1#2{.42#1}\epsffile{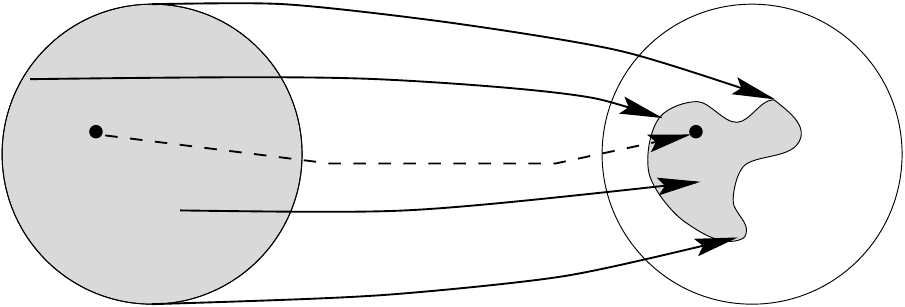}
\]
\end{Theorem}

(The term \emph{null-homotopic} that appears here refers to a map
that can be deformed to a constant map.) 

\proofheader{Proof of the equivalences} 
(Br1)$\Longrightarrow$(Br2) is trivial, since $S^{d-1}\sse B^d$.
   
For (Br2)$\Longrightarrow$(Br3) let
$h\:S^{d-1}\times[0,1]\to S^{d-1}$ be a null-homotopy for $f$, i.\,e.,
a continuous map that interpolates between our original map~$f$ and a
constant map, with 
$h(\xx,0)=f(\xx)$ and $h(\xx,1)=\xx_0$ for all $\xx\in S^{d-1}$. 
From this we
construct a continuous map $F\: B^d\to S^{d-1}$ that extends $f$, by
\[
F(\xx)\assg
\begin{cases}
h(\frac{\xx}{| \xx|},2-2| \xx|)    &\text{ if } \frac12\le|\xx|\le 1,\\
\xx_0                                &\text{ for }|\xx|\le\frac12. 
\end{cases}
\]
\[
\input{EPS/brouwer2-3a.pstex_t}
\]
This map is continuous, and by (Br2) it has a fixed point,
which must lie in the image, that is, in~$S^{d-1}$. 

For the converse, (Br3)$\Longrightarrow$(Br2),
let $f\: B^d\to S^{d-1}$ be continuous. Then
the restriction $f|_{S^{d-1}}$ is null-homotopic, since
$h(\xx;t)\assg f((1-t)\xx)$ provides a null-homotopy. 
Thus, by (Br3) the map $f|_{S^{d-1}}$ has a fixed point, hence so does~$f$. 

Finally, we get (Br2)$\Longrightarrow$(Br1): If $f\: B^d\to B^d$ 
has no fixed point, then we set 
$g(\xx)\assg\frac{f(\xx)-\xx}{| f(\xx)-\xx|}$.
This defines a map $g\:B^d\to S^{d-1}$
that has a fixed point $\xx_0\in S^{d-1}$ by (Br2),
with $\xx_0=\frac{f(\xx_0)-\xx_0}{| f(\xx_0)-\xx_0|}$.
But this implies $f(\xx_0) =\xx_0(1+t)$ for $t\assg{| f(\xx_0)-\xx_0|}>0$,
and this is impossible for $\xx_0\in S^{d-1}$. 
\endproof

In the following we use the unit cube $[0,1]^d$ in place of the
ball~$B^d$: It should be clear that the Brouwer fixed point theorem
equally applies to self-maps of any domain $D$ that is homeomorphic to the ball~$B^d$,
resp.\ of the boundary $\partial D$ of such a domain.

\proofheader{Proof of the Brouwer fixed point theorem 
        \rm(``HEX $\Longrightarrow$ (Br1)'')}  
If $f\: [0,1]^d\to [0,1]^d$ has no fixed point, then for some $\epsilon>0$ we have that
$| f(\xx)-\xx|_\infty\ge \epsilon$ for all $\xx\in [0,1]^d$ 
(namely, one can take 
$\epsilon\assg\min\{| f(\xx)-\xx|_\infty\sep \xx\in [0,1]^d\}$, 
which exists since $[0,1]^d$ is compact). 

Furthermore, any continuous function on the compact
set $[0,1]^d$ is uniformly continuous (see e.g.\ Munkres \cite[\S27]{Munkres1}), hence there exists some $\delta>0$
such that $|   \xx -  \xx' |_\infty <\delta$ 
implies   $| f(\xx)-f(\xx')|_\infty <\epsilon$. 
We take $\delta <\epsilon$ (without loss of generality), and then
choose $n$ with ${1\over n}<\delta$.

{From}~$f$, we now define a $d$-coloring of~$H(n,d)$, by setting
\[
\kappa(\vv)\assg 
{\textstyle\min\{i\sep |f_i(\frac{\vv}{n})-\frac{v_i}{n}|\geq\epsilon\}}
\]
for the interior vertices $\vv\in H(n,d)$, where $f_i$ denotes the $i$th component of~$f$.
This is well-defined, since $\frac{\vv}{n}\in [0,1]^d$, and thus
the absolute value of at
least one component of $f(\frac{\vv}{n})-\frac{\vv}{n}$ has to be at
least $\epsilon$. Now, the $d$-dimensional
HEX theorem guarantees a chain $\vv^0,\vv^1,\dots,\vv^N$ of
vertices of color~$i$, for some $i$, where $v^0_i=0$ and $v^N_i=n$. Furthermore, we know  that
$|f_i({\vv^k\over n})-{v^k_i\over n}|\geq\epsilon$ for $0\leq k\leq N$.
Also, at the ends of the chain we know the signs:
\begin{itemize}
\item[]
$f  (\tfrac{\vv^0}{n})\in [0,1]^d$ implies 
$f_i(\tfrac{\vv^0}{n})\ge0$ and hence
$f_i(\tfrac{\vv^0}{n}) - \tfrac{v^0_i}{n}\geq\epsilon$, and
\item[]
$f  (\tfrac{\vv^N}{n})\in [0,1]^d$ implies 
$f_i(\tfrac{\vv^N}{n})\le1$ and hence
$f_i(\tfrac{\vv^N}{n}) - \tfrac{v^N_i}{n}\leq -\epsilon$.
\end{itemize}
It follows that for some $k\in\{1,2,\dots,N\}$ we must have a sign change:
\begin{itemize}
\item[]
$f_i(\tfrac{\vv^{k-1}}{n}) - \tfrac{v^{k-1}_i}{n}\ge \epsilon$ and
$f_i(\tfrac{\vv^ k   }{n}) - \tfrac{v^k    _i}{n}\le-\epsilon$.
\end{itemize}
All these facts  taken together provide a contradiction, since
\begin{itemize}
\item[]
$|\tfrac{\vv^{k-1}}{n}-\tfrac{\vv^k}{n}|_\infty = \tfrac{1}{n}< \delta,$
\end{itemize}
whereas
\[
| f  (\tfrac{\vv^{k-1}}{n}) - f  (\tfrac{\vv^k}{n}) |_\infty \ge
  | f_i(\tfrac{\vv^{k-1}}{n}) - f_i(\tfrac{\vv^k}{n})  |        \ge
2\epsilon - |\tfrac{v^{k-1}_i}{n}-\tfrac{v^k_i}{n}|        \ge
2\epsilon - \tfrac{1}{n} > 2\epsilon - \delta >\epsilon.
\]
\endproof

\proofheader{Proof that the Brouwer fixed point theorem implies the HEX theorem
\rm(``Br1 $\Longrightarrow$ HEX'')}
Assume we have a coloring of~$H(n,d)$. We use it to define a map
$[0,n]^d\to [0,n]^d$, as follows: On the points in
$\{0,1,\dots,n\}^d$ we define
\[
f(\vv)=\begin{cases}
            \vv+\ee_i &\text{if $\vv$ has color $i$, 
            and there is a path on vertices of color $i$}\\[-4pt]
                  &\text{that connects $\vv$ to a vertex $\ww$ with $w_i=0$}\\
            \vv-\ee_i &\text{if $\vv$ has color $i$, but there is no such path.}
        \end{cases}
\]
If for the given coloring there is no winning path for HEX, then these definitions do not map any
point $\vv$ outside $[0,n]^d$. Hence this by linear extension
defines a simplicial map $f\: [0,n]^d\to [0,n]^d$ on the
simplices of the triangulation $\Delta(n,d)$ that we have considered before. 

The following two observations now give us a contradiction, showing that
this $f$ cannot have a fixed point:
\begin{compactitem}[ $\bullet$]
\item If $\Delta=\conv\{\vv^0,\vv^1,\vv^2,\dots,\vv^d\}\subseteq\R^d$ 
is a simplex and $f\: \Delta\to \R^d$ is a linear map defined by
$f(\vv^i)=\vv^i+\ww^i$, then $f$ has a fixed point on $\Delta$ if and
only if $\colzero\in\conv\{\ww^0,\dots,\ww^d\}$.
\item If $\vv,\vv'$ are adjacent vertices, then we cannot get
$f(\vv)=\vv-\ee_i$ and $f(\vv')=\vv'+\ee_i$.
Hence for each simplex of $\Delta(n,d)$, all the vectors $\ww^i$ lie
in one orthant of~$\R^d$!
\endproof
\end{compactitem}

\subsection{The joy of HEX: Who wins?}

So, who can win the $2$-dimensional
HEX game? A simple but ingenious argument due to John Nash,
known as “stealing a strategy,”
shows that on a square board the first player (``White'') always has
a winning strategy. In the following we first define winning
strategies, then show that one of the players has one, and finally
conclude that the first player has one. Still: The proof will be
non-constructive, and we don't know how to win HEX. So, the game
still remains interesting~\ldots

\begin{Definition}
A \emph{strategy} is a set of rules that tells one of the players which move to
choose (i.\,e., which tile to color) for every legal position on the board.
A \emph{winning strategy} here guarantees to lead to a win, 
starting from an empty board, for all
possible moves of the opponent.

A \emph{position} of the HEX game is a board on which some tiles may
have been colored white or black, together with the information who
moves next (unless all tiles are colored). A position is \emph{legal}
if it can occur in a HEX game: That is, if either White moves next,
and the numbers of white and black tiles agree, or if Black moves
next, and White has one more tile.

A \emph{winning position for White} is a position such that White has
a \emph{winning strategy} that tells him 
how to proceed (for arbitrary
moves of Black) and guarantees a win. Similarly, \emph{a winning
position for Black} has a \emph{winning strategy} that guarantees to
lead Black to a win.
\end{Definition}

\begin{Lemma}\label{nodraw}
Every (legal) position for HEX is either a winning
position for White or a winning position for Black.
\end{Lemma}

\proof
Here we proceed by induction on the number $g$ of ``grey'' tiles
(i.\,e., ``free'' positions on the board). If
no grey tiles are present $(g=0)$, then one of the players has won---by the HEX theorem.

If $g>0$ and White is to move, then any move that White could choose
reduces $g$, and thus (by induction) produces a winning position for
one of the players. If there is a move that leads to a winning
position for White, then this is really nice and great for White:
this makes the present position into a winning position for White,
and any such move can be used for a winning position for White.
Otherwise---too bad: If every possible move for White produces a
winning position for Black, then we are at a winning position for
Black already.

And the same argument applies for $g>0$ if Black is to move.
\endproof         

Of course, the argument given here is \emph{much} more general:
essentially we have proved that for any finite deterministic 2-person
game without a draw and with ``complete information'' there is a
winning strategy for one of the players. (This is a theorem of
Zermelo, which was rediscovered by von Neumann and Morgenstern). 
Furthermore, for games where a draw is possible either
one player has a winning strategy, or \emph{both} players can force a draw. We refer to 
Exercise~\ref{ex:drawplay}, and to Blackwell \& Girshick \cite[p. 21]{Blackwell}.

For HEX, Lemma~\ref{nodraw} shows that at the beginning (for the
starting position, where all tiles are grey, and White is to move),
there is a winning strategy either for White or for Black. But who is the winner?

Our first attempt might be to follow the proof of Lemma~\ref{nodraw}.
Only for the $2\times 2$ board this can be done:

\hspace{-10mm}\input{EPS/hexstrategy3corr.pstex_t}  

\noindent
In this drawing, you can decide for every position whether it
is a winning position for White or for Black, starting with the
bottom row ($g=0$) that has three winning positions for each player,
ending at the top node ($g=4$), which turns out to be a winning position for~White.

For larger boards, this approach is hopeless---after all, there are
$\binom{n^2}{\lfloor n^2/2\rfloor}$ final positions to classify for
``$g=0$,'' and from this one would have to work one's way up to the top node
of a huge tree (of height $n^2$). Nevertheless, people have worked out 
winning strategies for White on the $n\times n$ boards for $n\leq 5$
(see Gardner \cite{Gardner-hex}).

\begin{Theorem} 
For the HEX game played on a HEX board with equal side lengths,
White (the first player) has a winning strategy. 
\end{Theorem}

\proof
Assume not. Then by Lemma~\ref{nodraw} Black has a winning strategy. But then
White can start with an arbitrary move, and then---using the
symmetry of the board and of the rules---just ignore his first
tile, and follow Black's winning strategy ``for the second player.''
This strategy will tell White always which move to take. 
Here the ``extra'' white tiles cannot hurt White: If the
move for White asks to occupy a tile that is already white, 
then an arbitrary move is fine for White.
But this ``stealing a strategy'' argument produces a winning strategy for White, contradicting our
assumption!
\endproof

\begin{bibpar}
Gale's beautiful paper \cite{Gale-hex} was the source and
inspiration for our treatment of Brouwer's fixed point theorem in
terms of the HEX game. Nash's analysis for the winning strategies    
for HEX is from Gardner's classical account in \cite{Gardner-hex},
some of which reappears in Milnor's \cite{Milnor-nash}. 
See also the accounts in Jensen \& Toft \cite[Sect. 17.14]{Jensen-hex},
and in Berlekamp, Conway \& Guy \cite[p. 680]{BCG2},
where other cases of ``strategy stealing'' are discussed.
(A theoretical set-up for this is in 
Hales \& Jewett \cite[Sect. 3]{HalesJewett}.)

The traditional
combinatorial approach to the Brouwer fixed point theorem is via
Sperner's lemma \cite{Sperner}; see e.g.\ Exercise \ref{ex:Sperner} below and 
the presentation in \cite{AZ5}. Lov\'asz's \cite{Lov80} matroid version of Sperner's lemma
in Exercise \ref{ex:Lovasz-matroid} was further generalized by
Lindstr\"om \cite{Lin81}. Kry\'nski \cite{Kry90}, however, showed that these
results can easily be derived from earlier results.

A more geometric version of the
combinatorial lemmas is given by Mani \cite{Mani-lemma}.
\end{bibpar}

\begin{exs}
\item
Stir your coffee cup. Show that  
the (moving, but flat) surface has at every moment at least one point that 
stands still (has velocity zero). 
\item 
Prove that if you tear a sheet of paper from your notebook, crumble it
into a small ball, and put that down on your notebook,
then at least one point of the sheet comes to rest exactly on top of
its original position.\\
Could it happen that there are exactly two such points?
\item In the proof of the Brouwer fixed point theorem (Thm. \ref{t:brouwer}, 
(Br2)$\Longrightarrow$(Br3)), we could have tried to simply put
$F(\xx)\assg h(\frac{\xx}{| \xx|},1-| \xx|)$. Is this continuous? 
\item\label{ex:Sperner} (a) Prove ``Sperner's Lemma'' \cite{Sperner}: Let $\Delta$ be a triangulation of the 
$d$-dimensional sphere and let us color the vertices of~$\Delta$ using 
$d+1$ colors. Then $\Delta$ has an even number of colorful 
facets (meaning $d$-faces containing vertices of all colors). \\
(b) Show that Sperner's Lemma implies the Brouwer fixed point theorem.
\item\label{ex:Lovasz-matroid} (a) Let $\Delta$ be a triangulation of a 
$d$-dimensional manifold with vertex set $V$. 
Assume that a matroid $M$ of rank $d+1$ without loops
is defined on $V$. If $\Delta$ has a facet that is a basis of~$M$ then
it has at least two such facets. (Lov\'asz \cite{Lov80})\\ 
 (b) Show that part (a)  implies Sperner's Lemma, and hence also Brouwer's theorem. 
\item Let $B_E =2^{E}\setminus\{\emptyset, E\}$ be the poset of all proper subsets of a  
finite set $E$, ordered by containment.
Show that if an order-preserving map $f\: B_{E} \rightarrow B_{E}$
does not have a fixed point then it is 
surjective, and hence an automorphism.
\item Let $P=B_E \setminus\{A\}$, for  some proper subset $A$.\\
(a) Give a quick proof that $P$ has the \emph{fixed point property}, 
meaning that any order-preserving self-map has a fixed point.\\
(b) Give a slow proof, not using topology, that $P$ has the fixed point 
property. 
\item 
For HEX on a $3\times3$ board, how large is the tree of possible positions?
\item 
Can you write a computer program that plays HEX and wins (sometimes) \cite{Browne-hex_strat}?
\item 
For $d$-dimensional HEX, is there always some
``short'' winning path? Show that for every $d\geq 2$ there
is a constant $c_d$ such that for all $n$ there is a final
configuration such that only one player wins, but his shortest path
uses more than $c_d\cdot n^d$ tiles.
\item 
Construct an algorithm that, for given $\epsilon >0$ and
$f\: [0,1]^2\to [0,1]^2$, calculates a point $x_0\in [0,1]^2$ with
$| f(x_0)-x_0| <\epsilon$. \cite[p. 827]{Gale-hex} 
\item\label{ex:drawplay}
If in a complete information two player game a draw is possible,
argue why either one of the players has a winning strategy, or 
\emph{both} can force at least a draw.
\item\label{exer:jct}
Prove that for $2$-dimensional HEX, not both players can win!
For this, prove and use the ``polygonal Jordan curve theorem'':
any simple closed polygon in the plane uniquely divides the plane into an
``inside'' region and an ``outside'' region.      \\
(The general Jordan curve theorem for simple ``Jordan arcs'' in the
plane has extensive discussions in many books; see for example 
Munkres \cite{Munkres1}, Stillwell \cite[Sect. 0.3]{Stillwell},
or Thomassen~\cite{Thomassen-J}.)
\item 
On an $(m\times n)$-board that is not square 
(that is, $m\neq n$), the player who gets the longer sides, and hence the shorter
distance to bridge by a winning path, has a winning strategy. Our
figure illustrates the case of a $(6\times 5)$-board, where the claim
is that Black has a winning strategy.
\begin{itemize}
\item[(i)] Show that for this, it is sufficient to consider the case
where $m=n+1$ (i.\,e., the second player Black, who gets the longer
side, has a sure win).
\[
\def\epsfsize#1#2{.36#1}\epsffile{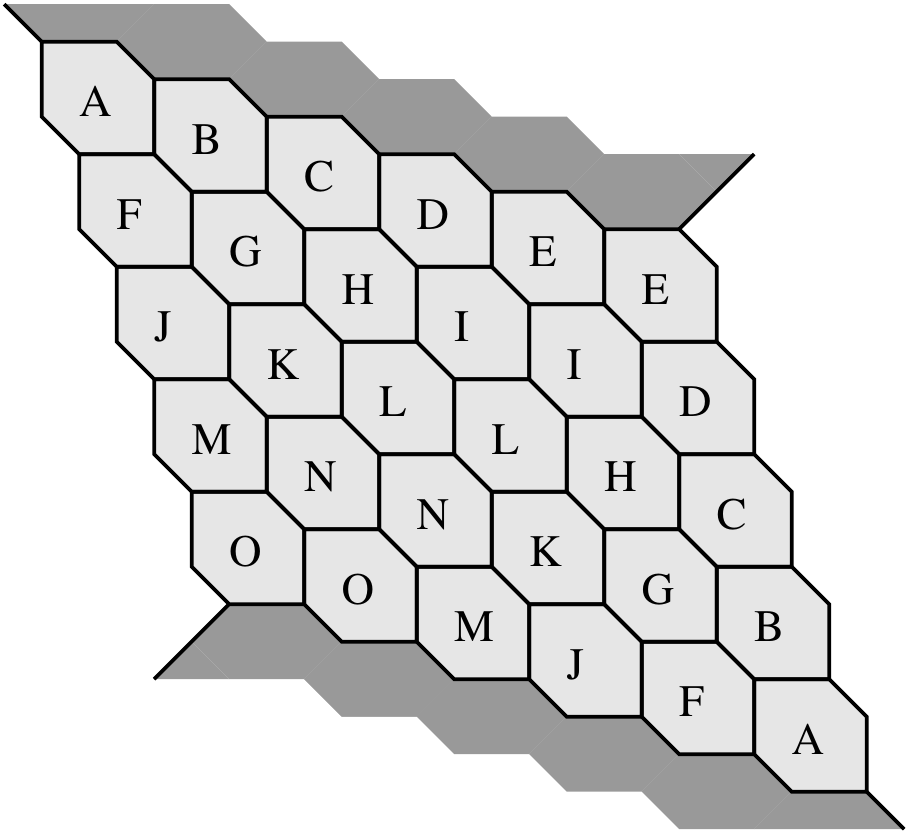}
\]
\item[(ii)] Show that in the situation of (i), Black has the
following winning strategy. Label the tiles in the ``symmetric''
way that is indicated by the figure, such that
there are two tiles of each label. The strategy
for Black is to always take the second tile that has the same label
as the one taken by White. 
Why will this strategy always win for Black?
(Hint: You will need the Jordan curve theorem.)\\
(This is in Gardner \cite{Gardner-hex}
and in Milnor \cite{Milnor-nash}, but neither source gives the proof.
You'll have to work yourself!) 
\end{itemize}
\end{exs}

\bigskip
\section{Piercing multiple intervals}\label{s:kais}
\subsection{Packing number and transversal number }

Let $\Sfrak$ be a system of subsets of a ground set~$X$; both $\Sfrak$
and $X$ may generally be infinite. The \emph{packing number%
\index{packing number}%
\index{number!packing}%
} of~$\Sfrak$, usually denoted by $\nu(\Sfrak)$ and often also
called the \emph{matching number%
\index{matching number!see{packing number}}%
\index{number!matching!see{packing number}}%
}, is the maximum cardinality of a system of pairwise disjoint
sets in $\Sfrak$:
\[
\nu(\Sfrak)=\sup\{|\MM|\sep \MM\subseteq\Sfrak,\, M_1\cap M_2=\es
\mbox{ for all $M_1, M_2\in \MM$, $M_1\neq M_2$}\}.%
\index{0nS@$\nu(\Sfrak)$}%
\]
\[
\includegraphics[height=22mm]{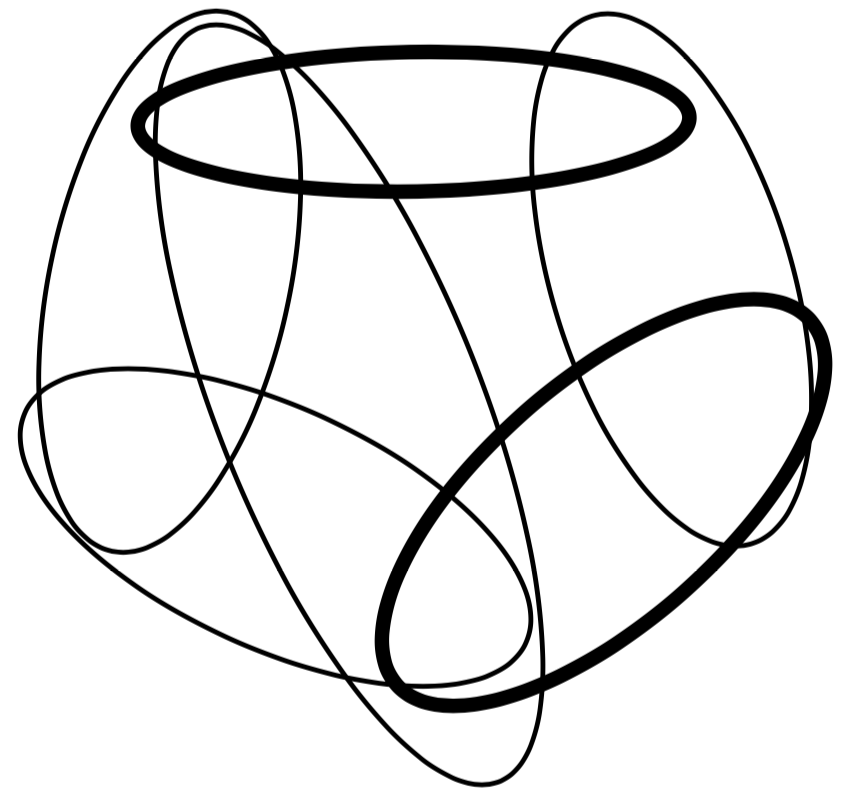} 
\]
The \emph{transversal number} or \emph{piercing number}
of~$\Sfrak$ is the smallest number of points of~$X$ that capture all
the sets in $\Sfrak$:
\[
\tau(\Sfrak)=\min\{|T|\sep T\subseteq X,\, S\cap T\neq\es\mbox{ for
all $S\in\Sfrak$}\}.
\index{0tS@$\tau(\Sfrak)$}%
\]
\[
\includegraphics[height=22mm]{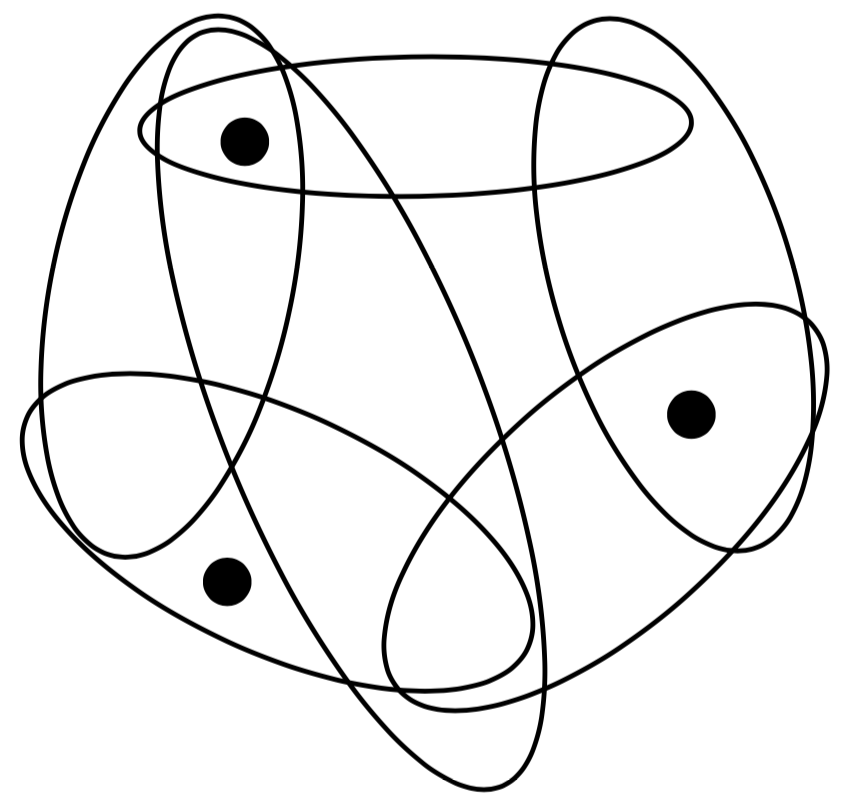}
\]
A subsystem $\MM\subseteq \Sfrak$ of pairwise disjoint sets is usually called a
\emph{matching} (this refers to the graph-theoretical matching,
which is a system of pairwise disjoint edges), and a set
$T\subseteq X$ intersecting all sets of~$\Sfrak$ is referred to as a
\emph{transversal} of~$\Sfrak$. Clearly, any transversal is at
least as large as any matching, and so always
\[
\nu(\Sfrak)\leq \tau(\Sfrak).
\]
In the reverse direction, very little can be said in general,
since $\tau(\Sfrak)$ can be arbitrarily large even if $\nu(\Sfrak)=1$.
As a simple geometric example, we can take the plane
as the ground set of~$\Sfrak$ and let the sets of~$\Sfrak$ be lines
in general position. Then $\nu=1$, since every two lines intersect,
but $\tau\geq \frac 12|\Sfrak|$, because no point is contained in
more than two of the lines.

One of the basic general questions in combinatorics asks for
interesting special classes of set systems where the transversal
number can be bounded in terms of the matching
number.\footnote{This kind of problem is certainly not
restricted to combinatorics. For example, if $\Sfrak$ is the system
of all open sets in a topological space, $\tau(\Sfrak)$ is the
minimum size of a dense set and is called the \emph{density%
\index{density}%
},
while $\nu(\Sfrak)$ is known as the \emph{Souslin number%
\index{Souslin number}%
\index{number!Souslin}%
}
or \emph{cellularity%
\index{cellularity}%
} of the space. In 1920, Souslin asked
whether a linearly ordered topological space exists (the
open sets are unions of open intervals) with countable $\nu$ but
uncountable $\tau$. It turned out in the 1970s that the answer 
depends on the axioms one is willing to assume beyond the
usual (ZFC) axioms of set theory. For example,
it is yes if one assumes the continuum hypothesis; see e.\,g.\
\cite{Engelking}.}
Many such examples come from geometry.
Here we restrict our attention
to one particular type of systems, the \emph{$d$-intervals},
where the best results have been obtained by topological methods.

\heading{Fractional packing and transversal numbers. }
Before introducing $d$-intervals, we mention another important
parameter of a set system, which always lies between $\nu$ and
$\tau$ and often provides useful estimates for $\nu$ or $\tau$.
This parameter can be introduced in two seemingly different ways.
For simplicity, we restrict ourselves to finite set systems
(on possibly infinite ground sets). A \emph{fractional packing%
\index{fractional packing}%
\index{packing!fractional}%
}
for a finite set system $\Sfrak$ on a ground set $X$
is a function $w\: \Sfrak\to [0,1]$
such that for each $x\in X$, we have $\sum_{S\in\Sfrak\sep x\in S} w(S)\leq 1$.
The \emph{size} of a fractional packing $w$ is 
$\sum_{S\in\Sfrak} w(S)$, and the \emph{fractional packing number}
$\nu^*(\Sfrak)%
\index{0n*S@$\nu^*(\Sfrak)$}%
$ is the supremum of the sizes of all fractional packings for~$\Sfrak$.
So in a fractional packing, we can take, say, one-third of one set
and two-thirds of another, but at each point, the fractions
for the sets containing that point must add up to at most~1.
We always have $\nu(\Sfrak)\leq \nu^*(\Sfrak)$, since a packing $\MM$ defines
a fractional packing $w$ by setting $w(S)=1$ for $S\in\MM$ and
$w(S)=0$ otherwise.

Similar to the fractional packing, one can also introduce
a fractional version of a transversal. A \emph{fractional transversal%
\index{fractional transversal}%
\index{transversal!fractional}%
}
for a (finite) set system $\Sfrak$ on a ground set $X$ is a function
$\varphi\:X\to [0,1]$ attaining only finitely many nonzero values
such that for each $S\in \Sfrak$, we have $\sum_{x\in S} \varphi(x)\geq 1$.
The size of a fractional transversal $\varphi$ is
$\sum_{x\in X}\varphi(x)$, and the \emph{fractional transversal number}
$\tau^*(\Sfrak)%
\index{0t*S@$\tau^*(\Sfrak)$}%
$ is
the infimum of the sizes of fractional transversals.

By the duality theorem of linear programming (or by the theorem about separation
of disjoint convex sets by a hyperplane), it follows that $\nu^*(\Sfrak)=\tau^*(\Sfrak)$
and thus that
\[
	\nu(\Sfrak) \ \le\ \nu^*(\Sfrak)\ =\ \tau^*(\Sfrak)\ \le\ \tau(\Sfrak)
\]
for any finite set system~$\Sfrak$.

When trying to bound $\tau$ in terms of~$\nu$, in many instances
it proved very useful to bound $\nu^*$ as a function of~$\nu$
first, and then $\tau$ in terms of~$\tau^*$. The proof presented
below follows a somewhat similar approach.

\subsection{The \emph{d}-intervals}

Let $I_1,I_2,\dots,I_d$ be disjoint parallel segments in the plane. 
(We may assume without loss of generality that they are horizontal
unit length intervals at distinct heights/$y$-coordinates.)
A set $J\subset \bigcup_{i=1}^d I_i$ is a 
\emph{$d$-interval\index{interval, $d$-interval}\index{dinterval@$d$-interval}}
if it intersects each $I_i$ in a closed interval.
We denote this intersection by $J_i$ and call it the \emph{$i$th component} of~$J$. 
The following drawing shows a $3$-interval:
\[
\includegraphics[height=22mm]{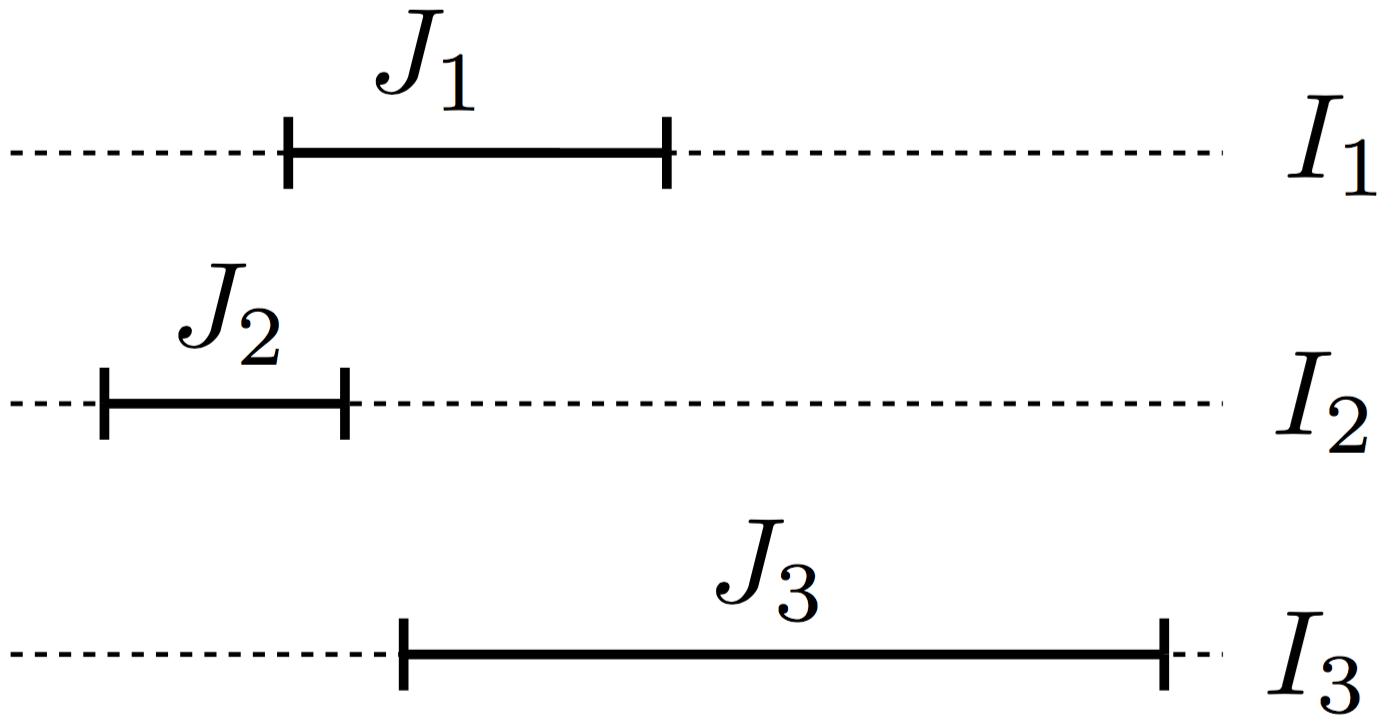}
\]
Intersection and piercing for $d$-intervals are taken in the
set-theoretical sense: Two $d$-intervals intersect if,
for some $i$, their $i$th components intersect.

The $1$-intervals, which are just intervals in the usual sense,
behave nicely with respect to
packing and piercing, as for any family $\FF$ of intervals, we
have $\nu(\FF)=\tau(\FF)$. (This is well-known and easy to prove: Exercise \ref{exerc:intervals}!))
This, however, does not extend to $d$-intervals.
For example, the family $\FF$ of three $2$-intervals
\[
\includegraphics[height=25mm]{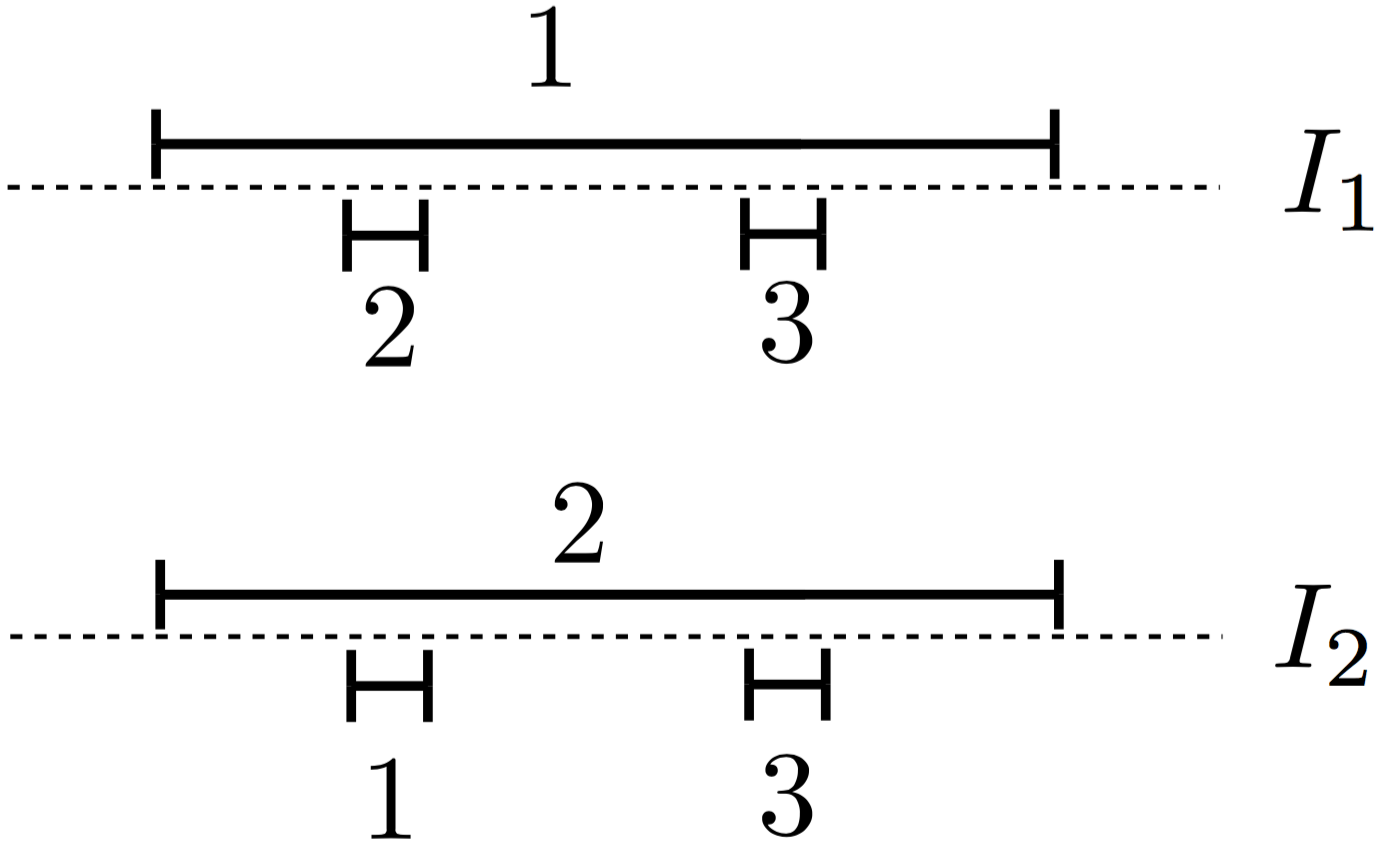}
\]
has $\nu(\FF)=1$ while $\tau(\FF)=2$. By taking multiple
copies of this family, one obtains families with $\tau=2\nu$
for all values of~$\nu$.

Gy\'arf\'as \& Lehel \cite{GyarfasLehel0}
showed by elementary methods that for any $d$ and any family
$\FF$ of $d$-intervals, $\tau(\FF)$ can be bounded by a function
of~$\nu(\FF)$ (also see \cite{GyarfasLehel}). 
Their function was rather large (about
$\nu^{d!}$ for $d$ fixed). After an initial breakthrough
by Tardos \cite{Tardos-2int}, who proved $\tau(\FF)\leq
2\nu(\FF)$ for any family of $2$-intervals, Kaiser
\cite{Kaiser} obtained the following result:

\begin{Theorem}[The Tardos--Kaiser theorem on $d$-intervals]
\label{t:d-int}%
\index{Tardos--Kaiser theorem|thmref{t:d-int}}%
\index{theorem!Tardos--Kaiser|thmref{t:d-int}}%
Every family $\FF$ of {$d$-in}\-ter\-vals,
$d\geq 2$, has a transversal of
size at most $(d^2-d)\cdot \nu(\FF)$.
\end{Theorem}
 
Here we present a proof using the Brouwer fixed point theorem.
Alon \cite{Alon-piercing} found a short non-topological proof
of the slightly weaker bound $\tau(\FF)\leq 2d^2\nu(\FF)$.

\proof Let $\FF$ be a fixed system of $d$-intervals
with $\nu(\FF)=k$,
and let $t=t(d,k)$ be a suitable (yet undetermined) integer.
The general plan of the proof is this: Assuming that there
is no transversal of~$\FF$ of size $dt$, we show by a topological method
that the fractional packing number $\nu^*(\FF)$ is at least $t+1$.
Then a simple combinatorial argument proves that
the packing number $\nu(\FF)$ is at least $\frac {t+1}d$,
which leads to $t<d^2\cdot\nu(\FF)$. A sharper combinatorial
reasoning in this step leads to the slightly better bound
in the theorem.

Our candidates for a transversal of~$\FF$ are all
sets $T$ with each $T_i=T\cap I_i$ having exactly $t$ points;
so $|T|=td$. For technical reasons, we also permit
that some of the $t$ points in $I_i$ coincide, so $T$ can be a
multiset.

The letter $T$ could also abbreviate a \emph{trap}.
The trap is
set to catch all the $d$-intervals in~$\FF$, but if it is not
set well enough, some of the $d$-intervals can escape.
Each of them escapes through a hole in the trap, namely
through a \emph{$d$-hole}. The points of~$T_i$
cut the segment $I_i$ into $t+1$ open intervals
(some of them may be empty), and these
are the \emph{holes in $I_i$}; they are numbered 1 through
$t+1$ from left to right.
A $d$-hole consists of $d$~holes, one in each $I_i$. The \emph{type}
of a $d$-hole $H$ is the set
 $\{(1,j_1),(2,j_2),\dots,(d,j_d)\}$, where
$j_i\in [t{+}1]$ is the number of the hole
in $I_i$ contained in~$H$.
A $d$-interval $J\in\FF$ \emph{escapes}
through a $d$-hole $H$ if it is contained in the union of its
holes. The drawing shows a $3$-hole, of type
$\{(1,2),(2,4),(3,4)\}$, and a $3$-interval
escaping through it:
\[
\includegraphics[height=25mm]{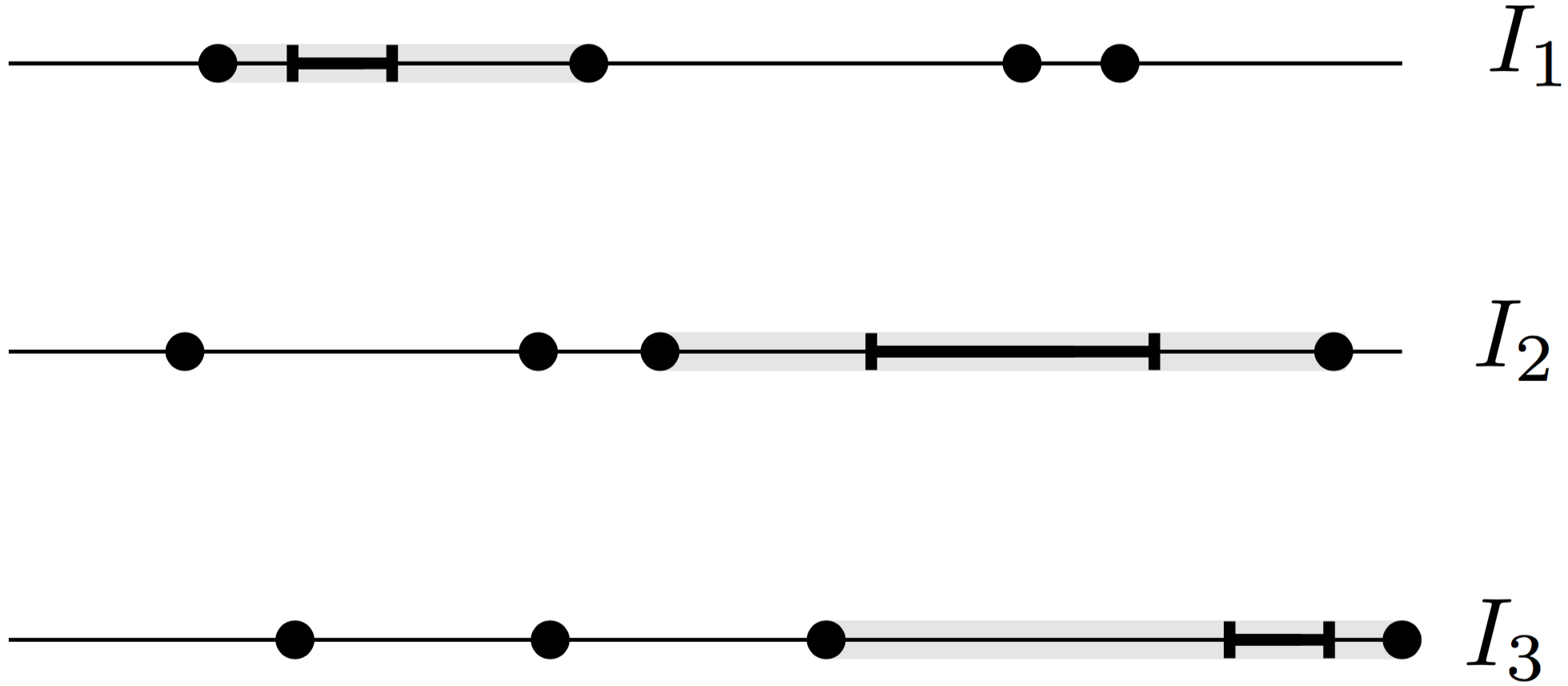}
\]
Let $\HH_0$ be the hypergraph with vertex set $[d]\times
[t{+}1]$ and with edges being all possible types of $d$-holes;
for example, the hole in the picture yields the edge
$\{(1,2),(2,4),(3,4)\}$.
So $\HH_0$ is a complete $d$-partite $d$-uniform
hypergraph.  
By saying that a $J\in\FF$ escapes through
an edge $H$ of~$\HH_0$, we mean that $J$ escapes through
the $d$-hole (uniquely) corresponding to~$H$.

Next, we define weights on the edges of~$\HH_0$;
these weights depend on the set~$T$ (and also on $\FF$, but this
is considered fixed).
The weight of an edge $H\in\HH_0$ is
\[
q_H=\sup\{\dist(J,T)\sep J\in\FF,\mbox{ $J$ escapes through $H$}\}.
\]
Here $\dist(J,T):=\min_{1\le i\le d}\{\dist(J_i,T_i)\}$ and
$\dist(J_i,T_i)$ is the distance of the $i$th component
of~$J$ to the closest point of~$T_i$.
Thus $q_H$ can be interpreted as the largest margin by which
some $d$-interval from $\FF$ escapes through~$H$. 
If no members of~$\FF$ escape through~$H$, we define
$q_H$ as~0. Note that this is the only case where
$q_H=0$. Otherwise, if anything escapes, it does so by a positive
margin, since we are dealing with closed intervals.

From the edge weights, we derive weights of vertices:
The weight $w_v$ of a vertex $v=(i,j)$
is the sum of the weights of the edges of~$\HH_0$
containing~$v$. These weights, too, are functions of~$T$;
to emphasize this, we write $w_v=w_v(T)$.

\begin{Lemma}\label{l:samewt}
For any $d\geq 1$, $t\geq 1$, and any $\FF$, there is
a choice of\ \,$T$ such that all the vertex weights
$w_v(T)$, $v\in [d]\times[t{+}1]$,
coincide.
\end{Lemma}

It is this lemma whose proof is topological. We postpone
that proof and finish the combinatorial part first.

Let us suppose that a trap $T$ was chosen as in the lemma,
with $w_v(T)= W$ for all~$v$. If $W=0$ then
$T$ is a transversal, since
all edge weights are 0 and  no $J\in\FF$ escapes.
So suppose that $W>0$.

Let $\HH=\HH(T)\subseteq \HH_0$,
the \emph{escape hypergraph} of~$T$, consist of
the edges of~$\HH_0$ with nonzero weights.
Note that
\begin{equation}\label{e:nuHnuF}
\nu(\HH)\leq \nu(\FF).
\end{equation}  
Indeed, given a matching
$\MM$ in $\HH$, for each edge $H\in\MM$ choose
a $J\in\FF$ escaping through $H$---this gives a matching in~$\FF$.

We note that the re-normalized edge
weights $\tilde q_H=\frac 1W\,q_H$ determine
a fractional packing in $\HH$ (since the weights at each
vertex sum up to~1). For the size of this fractional
packing, which is the total weight of all vertices, we find
by double counting
\[ 
\sum_{H\in\HH}\tilde q_H =
\frac 1d \sum_{H\in\HH}\sum_{v\in H} \tilde q_H =
\frac 1d \sum_{v\in [d]\times [t{+}1]}\frac{w_v}W=
\frac 1d \sum_{v} 1 = t+1.
\] 
As 	$\nu^*(\HH)$ is the supremum of the weights of all fractional packings,
and $\tilde q_H$ is a particular fractional packing, this yields
$\nu^*(\HH)\geq
\sum_{H\in\HH}\tilde q_H =t+1$.	
	
The last step is to show that $\nu(\HH)$ cannot be small
if $\nu^*(\HH)$ is large. Here is a simple argument
leading to a slightly suboptimal bound, namely
$\nu(\HH)\geq \frac 1d\,\nu^*(\HH)$.

Given a fractional matching $\tilde q$ of size $t+1$ in $\HH$, a matching
can be obtained by the following greedy procedure:
Pick an edge $H_1$ and discard all edges intersecting it,
pick $H_2$ among the remaining edges, etc., until all edges are
exhausted. The $\tilde q$-weight of~$H_i$ plus all the edges
discarded with it is at most $d=|H_i|$, 
while all edges together
have weight $t+1$. Thus, the number of steps, and also the size
of the matching $\{H_1,H_2,\dots\}$, is at least
$\lceil\frac {t+1} d\rceil$.

If we set $t=d\cdot\nu(\FF)$,
we get  $\nu(\HH)>\nu(\FF)$, which contradicts (\ref{e:nuHnuF}).
Therefore, for this choice of~$t$, all the vertex weights must be 0,
and $T$ as in Lemma~\ref{l:samewt} is a transversal
of~$\FF$ of size at most $d^2\nu(\FF)$.

The improved bound $\tau(\FF)\leq (d^2-d)\cdot\nu(\FF)$ for $d\ge 3$ follows
similarly using a theorem of F\"uredi \cite{Furedi-nunustar},
which implies that
the matching number of any $d$-uniform $d$-partite hypergraph
$\HH$ satisfies $\nu^*(\HH)\leq (d-1)\nu(\HH)$.
(For $d=2$, a separate argument needs to be used, based on
a theoreom of Lov\'asz stating that $\nu^*(G)\le\frac32\nu(G)$ for all graphs
$G$.)
The Tardos--Kaiser theorem~\ref{t:d-int} is proved.
\proofend

\heading{Proof of Lemma~\ref{l:samewt}. } Let $\sigma^t$ denote
the standard $t$-dimensional simplex in~$\R^{t+1}$, i.e.
the set $\{\xx\in\R^{t+1}\sep x_j\geq 0,\,x_1+\cdots+x_{t+1}=1\}$.
A point $\xx\in\sigma^t$ defines a $t$-point multiset
$\{z_1,z_2,\dots,z_t\}\subset [0,1]$, $z_1\leq z_2\leq\cdots\leq
z_t$, by setting $z_k=\sum_{j=1}^k x_j$. Here is a picture for
$t=2$:
\[
\includegraphics[height=40mm]{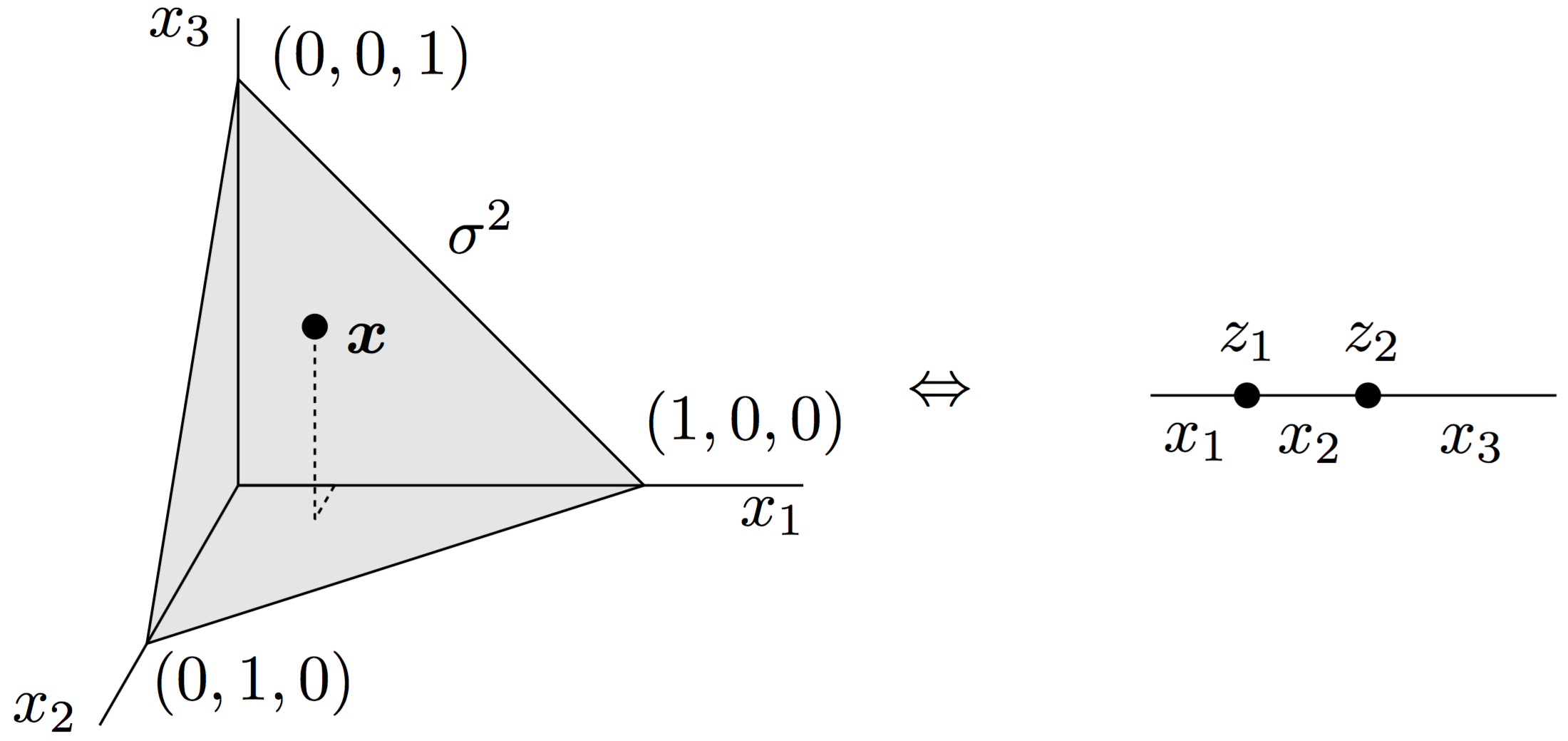}
\]
A candidate transversal $T$ with $t$ points in each $I_i$
can thus be defined
by an ordered $d$-tuple $(\xx_1,\dots,\xx_d)$ of points,
$\xx_i\in\sigma^t$, where $\xx_i$ determines $T_i$.
Such an ordered $d$-tuple can be regarded as a single point
 $\xx$ in the Cartesian product
$P=\sigma^t\times\sigma^t\times\cdots\times\sigma^t=
(\sigma^t)^d$. To each  $\xx\in P$, we have thus assigned
a candidate transversal $T(\xx)$.

For each vertex $v=(i,j)$ of the hypergraph $\HH_0$,
we define the function $g_{ij}\:P\to\R$ by
$g_{ij}(\xx)=w_{(i,j)}(T(\xx))$, where $w_v(T)$ is the
vertex weight. This is a \emph{continuous} function
of~$\xx$, 
since the edge weights $q_H$
and hence the vertex weights $w_{(i,j)}(T(\xx))$
change continuously when $T(\xx)$ 
moves---even if by this move new 
edges from $\FF$ escape, or fail to escape, through a hole:
If this is due to a small change of~$T(\xx)$,
then they escape, or fail to escape, by a narrow margin.

We note that for each $\xx$, the sum
\[
S_i(\xx)=\sum_{j=1}^{t+1} g_{ij}(\xx)
\]
is independent of~$i$; this is because $S_i(\xx)$ equals
the sum of the weights of all edges.
So we can write just $S(\xx)$ instead of~$S_i(\xx)$.

If there is an $\xx\in P$ with $S(\xx)=0$, then
all the vertex weights $ w_{(i,j)}(T(\xx))$ are 0
and we are done. Otherwise, we define the
normalized functions
\[
f_{ij}(\xx)= \frac 1{S(\xx)}\, g_{ij}(\xx).
\]
For each $i$, $f_{i1}(\xx),\dots,f_{i(t+1)}(\xx)$ are
nonnegative and sum up to 1, and so they are the coordinates of a
point in the standard simplex $\sigma^t$. All the maps
$f_{ij}$ together can be regarded as a map $f\:P\to P$.
To prove the lemma, we need to show that the image of~$f$
contains the point of~$P$ with all the $d(t+1)$
coordinates equal to $\frac 1{t+1}$.

The product $P$ is a convex polytope, and its nonempty faces
are exactly all Cartesian products
$F_1\times F_2\times\cdots\times F_d$, where the
$F_1,\dots,F_d$ are nonempty faces of the factors $\sigma^t,\dots,\sigma^t$
of $P=\sigma^t\times\sigma^t\times\cdots\times\sigma^t$
(Exercise~\ref{ex:polytprod}). We note that for any
face $F$ of~$P$, we have $f(F)\subseteq F$: Indeed,
any face $G$ of~$\sigma^t$ has the form
$G=\{\xx\in \sigma^t\sep x_i=0\mbox{ for all }i\in I\}$,
for some index set $I$, and the faces of~$P$ are products of
faces $G$ of this form.
 So it suffices to know that $f_{ij}(\xx)=0$ whenever
$(\xx_i)_j=0$. This holds, since $(\xx_i)_j=0$ means that
the $j$th hole in $I_i$ is empty, so nothing can escape
through that hole, and thus $f_{ij}(\xx)=0$.
The proof of Lemma~\ref{l:samewt} is now reduced to the following
statement.

\begin{Lemma} \label{l:polyt-surj}
Let $P$ be a convex polytope and let $f\:P\to P$ be
a continuous mapping satisfying $f(F)\subseteq F$ for each
face\footnote{In fact, it suffices to require $f(F)\subseteq F$
for each facet of~$P$ (that is, for each face of dimension $\dim(P)-1$), 
since each face is the intersection of some facets.}
$F$ of~$P$. Then $f$ is surjective.
\end{Lemma}

\proof 
Since the condition is hereditary for faces,
it suffices to show that each point $\yy$ in the interior of~$P$
has a preimage. For contradiction, suppose that some
$\yy\in\inter P$ is not in the image of~$f$.
For $\xx\in P$, consider the ray
that starts at $f(\xx)$ and passes through
$\yy$, and let $g(\xx)$ be the unique intersection of that
ray with the boundary of~$P$.
\[
\includegraphics[height=28mm]{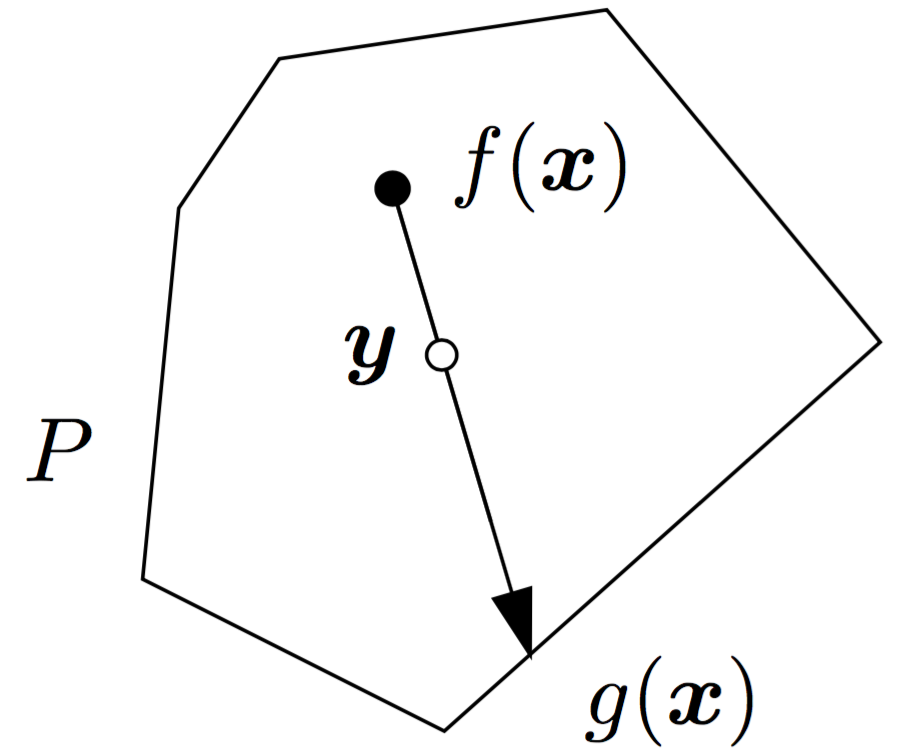}
\]
This $g$ is a well-defined and
continuous map $P\to P$, and by Brouwer's fixed point theorem,
there is an $\xx_0\in P$ with $g(\xx_0)=\xx_0$. The point $\xx_0$
lies on the boundary of~$P$, in some proper face $F$. But
$f(\xx_0)$ cannot lie in $F$, because the segment $\xx_0 f(\xx_0)$ passes
through the point $\yy$ outside~$F$---a contradiction.
\proofendSkip

\subsection{Lower bounds }
It turns out that the bound in Theorem~\ref{t:d-int} is
not far from being the best possible. In particular,
for $\nu(\FF)=1$ and $d$ large, the transversal number
can be near-quadratic in $d$, which is rather surprising.
For all $k$ and $d$, systems $\FF$ of $d$-intervals can be constructed
with $\nu(\FF)=k$ and
\[
\tau(\FF)\geq c\,{d^2\over(\log d)^2}\,k
\]
for a suitable constant $c>0$
(Matou\v{s}ek \cite{Matousek-d-int}). The construction
involves an extension of a construction
due to Sgall \cite{Sgall} of certain systems of set pairs.
Here we outline a (non-topological!)
proof of a somewhat simpler result concerning families of
\emph{homogeneous%
\index{interval, $d$-interval!homogeneous}%
\index{dinterval@$d$-interval!homogeneous}%
\index{homogeneous $d$-interval}%
} $d$-intervals,
which are unions of at most $d$ closed intervals on the real
line. These are more general than the $d$-intervals,
but an upper bound only slightly weaker than
Theorem~\ref{t:d-int} can be proved for them along the same lines
(Exercise~\ref{ex:homog}): $\tau\leq (d^2-d+1)\nu$.

\begin{Proposition}\label{p:d-int-lb}
	There is a constant $c>0$ such that 
for every $d\ge 2$ and $k\geq 1$, there exists a system $\FF$ of 
homogeneous $d$-intervals
with $\nu(\FF)=k$ and
\[
\tau(\FF)\geq c\,{d^2\over \log d}\,k.
\]
\end{Proposition}

\proof
Given $d$ and $k$, we want to construct a system $\FF$ of
homogeneous $d$-intervals. Clearly, it suffices to consider the
case $k=1$, since for larger $k$, we can take $k$ disjoint copies
of the $\FF$ constructed for $k=1$. Thus, we want an $\FF$ in
which every two $d$-intervals intersect and with $\tau(\FF)$
large.

In the construction, we will use homogeneous $d$-intervals
of a quite special form: Each component is either a single point
or a unit-length interval. First, it is instructive to see why 
we cannot get a good example if all the components are only
points. In that case, the family $\FF$ is simply a $d$-uniform
hypergraph (whose vertices happen to be points of the real line).
We require that any two edges intersect, and thus any edge is a transversal
and we have $\tau(\FF)\leq d$. 

For the actual construction, let $n$ and $N$ be integer parameters
(whose value will be set later). Let $V=[n]$
be an index set, and  $I_v$, for $v\in V$,
be auxiliary pairwise disjoint
unit intervals on the real line. In each $I_v$, we choose
$N$ distinct points $x_{v,i}$, $i=1,2,\dots,N$.

The constructed system $\FF$ will consist of homogeneous
$d$-intervals $J^1,J^2,\dots,J^N$. For each $i=1,2,\dots,N$,
we choose auxiliary sets $\emptyset\subset B_i\subseteq A_i\subseteq V$ and
then construct $J^i$ as follows: 
\[
J^i=\Big(\,\bigcup_{v\in B_i} I_v\Big) \cup
\{ x_{u,i}\sep u\in A_i\setminus B_i\}.
\]
The picture shows an example of~$J^1$ for $n=6$, $A_1=\{1,2,4,5\}$
and $B_1=\{2,4\}$:
\[
\def\epsfsize#1#2{1#1}\epsffile{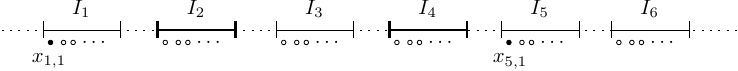}
\]
 The heart of the proof is the construction
of suitable sets $A_i$ and $B_i$ on the ground set~$V$. Since the $J^i$ should be
homogeneous $d$-intervals, we obviously require
\begin{enumerate}
\item[{\rm (C1)}]
For all $i=1,2,\dots,N$,
$\es\subset B_i\subseteq A_i$ and $|A_i|\leq d$.
\end{enumerate}
The condition that every two members of~$\FF$ intersect is implied by
the following: 
\begin{enumerate}
\item[{\rm (C2)}]
For all $i_1,i_2$, $1\leq i_1<i_2\leq N$, we have
$A_{i_1}\cap B_{i_2}\neq\emptyset$ or $A_{i_2}\cap B_{i_1}\neq\emptyset$
(or both).
\end{enumerate}
Finally, we want $\FF$ to have no small transversal. Since no
two $d$-intervals of~$\FF$ have a point component in common,
a transversal of size $t$ intersects no more than $t$ members
of~$\FF$ in their point components, and all the other members
of~$\FF$ must be intersected in their interval components.
Therefore, the transversal condition translates to
\begin{enumerate}
\item[{\rm (C3)}]
 Put $t=cd^2/\log d$ for a sufficiently small
constant $c>0$, and let $\BB=\{B_1,B_2,\dots,B_N\}$.
 Then  $\tau(\BB)\geq 2t$, and consequently $\tau(\BB')\geq t$
for any~$\BB'$ arising from $\BB$ by removing at most $t$ sets.
\end{enumerate}

A construction of sets $A_1,\dots,A_N$ and $B_1,\dots,B_N$
as above was provided by Sgall
\cite{Sgall}. His results give the following:

\begin{Proposition}\label{p:sgall}
Let $b$ be a given integer,  let $n\leq cb^2/\log b$
for a sufficiently small constant $c>0$, and let $B_1,B_2,\dots,B_N$
be $b$-element subsets of~$V=[n]$. Then there exist sets
$A_1,A_2,\dots,A_N$, with $B_i\subseteq A_i$, $|A_i|\leq 3b$,
and such that {\rm(C2)} is satisfied.
\end{Proposition}

With this proposition, the proof of Proposition~\ref{p:d-int-lb}
is easily finished. We set $b=\lfloor \frac d 3\rfloor$, 
$n=cb^2/\log b$, and we
let $B_1,B_2,\dots,B_N$ be all the $N=\binom{n}{b}$ subsets of
$V$ of size~$b$. We have $\tau(\{B_1,\dots,B_n\})=n-b+1$
and condition (C3) holds.
It remains to construct the sets $A_i$ according to 
Proposition~\ref{p:sgall}; then (C1) and (C2) are satisfied too.
The proof of Proposition~\ref{p:d-int-lb} is concluded
by  passing from the $A_i$ and $B_i$ to the system $\FF$
of homogeneous $d$-intervals as was described above.
\proofend

\proofheader{Sketch of proof of Proposition~\ref{p:sgall}}
Let $G=(V,E)$ be a graph on $n$ vertices of maximum degree $b$
with the following expander-type
property: For any two disjoint $b$-element subsets $A,B\subseteq
V$, there is at least one edge $e\in E$ connecting a vertex of
$A$ to a vertex of~$B$. (The existence of such a graph can be
easily shown by the probabilistic method; the constant
$c$ arises in this argument. See \cite{Sgall} for references.) 

For each $i$, let $v_i$ be an (arbitrary) element of
the set $B_i$, and let
\[
A_i= B_i\cup N(v_i)\cup \Big(V\setminus \bigcup_{u\in B_i}N(u)\Big),
\]
where $N(v)$ denotes the set of neighbors in $G$ of a vertex
$v\in V$. It is easy to check that $|A_i|\leq 3b$, and some
thought reveals that the condition (C2) is satisfied.
\proofend

\subsection{A Helly-type problem for \emph{d}-intervals}

Kaiser \& Rabinovich \cite{KaiserRabinovich} investigated
conditions on a family $\FF$ of $d$-intervals guaranteeing that
$\FF$ can be pierced by a ``multipoint,'' that is, $\tau(\FF)\le d$
and there is a transversal using one point of each~$I_i$.
They proved the following.

\begin{Theorem}[The Kaiser--Rabinovich theorem on $d$-intervals]\label{t:multipoint}
 Let $k=\lceil\log_2(d+2)\rceil$ and let
$\FF$ be a family of $d$-intervals such that any $k$ or fewer
members of~$\FF$ have a common point. Then $\FF$ can be pierced
by a multipoint.
\end{Theorem} 

Let's put this result into context:
The proof of the Kaiser--Tardos Theorem~\ref{t:d-int} sets out to show that there exists a 
transversal consisting of exactly $t$ points in each of the intervals~$I_i$, for a suitable~$t$. 
We eventually get that if every two $d$-intervals meet (that is, $\nu(\FF)=1$), then we can take $t<d$. 
The Kaiser--Rabinovich theorem says that if every $\lceil\log_2(d + 2)\rceil$
meet then $t<2$ suffices. The upcoming proof of Theorem \ref{t:multipoint} can be extended to
yield an interpolation between this result and the Kaiser--Tardos theorem:
If every $\lceil\log_b(d + 2)\rceil$ edges meet, then we can take $t<b$. 
For $b = d$ this yields the result of Kaiser--Tardos for $\nu(\FF)=1$.

\proof  
We use notation from the proof of
Theorem~\ref{t:d-int}.
We apply Lemma~\ref{l:samewt} with $t=1$, obtaining
a set $T$ with one point in each $T_i$ such that all the
$2d$ vertices of the escape hypergraph $\HH=\HH(T)$ have the same
weight~$W$. If $W=0$ we are done, so let us assume $W>0$.

By the assumption on $\FF$,
every $k$ edges of~$\HH$ share a common vertex.
We will prove the following claim for every $\ell$:
\begin{quote}
\emph{If every $\ell+1$ edges
of~$\HH$ have at least $m$ common vertices, then every
$\ell$ edges of~$\HH$ have at least $2m+1$ common vertices.}
\end{quote}
For $\ell=k$, the assumption holds with $m=1$, and so
by $(k-1)$-fold application of this claim, we get
that every edge of~$\HH$ ``intersects itself'' in
at least $2^k-1$ vertices, i.e. $d>2^k-2$. The claim thus
implies the theorem.

The claim is proved by contradiction. Suppose that
$\AA\subseteq\HH$ is a set of~$\ell$ edges such that
$C=\bigcap\AA$ has at most $2m$ vertices,
and let $\bar C :=\{(i,3-j)\sep (i,j)\in C\}$.
No edge $H\in\HH$ contains both $(i,1)$ and $(i,2)$,
thus also $C$ does not contain both $(i,1)$ and $(i,2)$,
and thus $\bar C$ is a subset of the complement of~$C$;
it is matched to $C$ by 
$(i,3-j) \leftrightarrow (i,j)$, and thus $|C| = |\bar{C}|$. 
	  
By the assumption, $\AA$
plus any other edge together intersect in at least $m$
vertices. Thus, any $H\in \HH\setminus \AA$ contains
at least $m$ vertices of~$C$, and consequently no more than $m$
vertices of~$\bar C$.

Let $U$ be the total weight of the vertices in $C$, and $\bar U$
the total weight of the vertices in $\bar C$. The edges in $\AA$
contribute solely to $U$, while any other edge $H$ contributes at
least as much to $U$ as to $\bar U$, and so $U>\bar U$.
But this is impossible since all vertex weights are identical and
$|C|=|\bar C|$. The claim, and Theorem~\ref{t:multipoint} too,
are proved.
\proofendSkip

An interesting open problem is whether
$k=\lceil\log_2(d+2)\rceil$ in Theorem~\ref{t:multipoint}
could be replaced by $k=k_0$ for some constant $k_0$ independent
of~$d$. The best known lower bound is $k_0\geq 3$. \\

\begin{bibpar}
Tardos \cite{Tardos-2int} proved the optimal bound
$\tau\leq 2\nu$ for 2-intervals by a topological argument
using the homology of suitable simplicial complexes.
Kaiser's argument \cite{Kaiser} is similar to the presented one,
but he proves Lemma~\ref{l:samewt} using a rather advanced
Borsuk--Ulam-type
theorem of Ramos \cite{Ramos} concerning continuous maps defined
on products of spheres. The method with Brouwer's theorem was used
by Kaiser \& Rabinovich \cite{KaiserRabinovich} for a proof of
Theorem~\ref{t:multipoint}.

Lemma \ref{l:polyt-surj} seems to be new in the version that
we give here, but it relates to a vast literature of 
``KKM-type lemmas,'' which starts with a paper by 
Knaster, Kuratowski, and Mazurkiewicz \cite{KKM29}
from 1929. We refer to Bárány \& Grinberg \cite{BaranyGrinberg} 
and the references given there, such as 
\url{http://mathoverflow.net/questions/67318}.

Alon's short proof \cite{Alon-piercing} of the bound
$\tau\leq 2d^2\nu$ for families of $d$-intervals
applies a powerful technique developed in Alon \&
Kleitman \cite{ak-pcs-92a}. For the so-called Hadwiger--Debrunner
$(p,q)$-problem solved in the latter paper, the quantitative
bounds are probably quite far from the truth. It would be
interesting to find an alternative topological approach to that problem,
which could perhaps lead to better bounds.
See, for example, Hell \cite{Hell-frac-Helly}.

The variant of the piercing problem for families of homogeneous%
\index{interval, $d$-interval!homogeneous}%
\index{dinterval@$d$-interval!homogeneous}%
\index{homogeneous $d$-interval}
 $d$-intervals has been considered simultaneously with $d$-intervals; see
\cite{GyarfasLehel} \cite{Tardos-2int} \cite{Kaiser} \cite{Alon-piercing}.
The upper bounds obtained for the homogeneous case
 are slightly worse:
$\tau\leq 3\nu$ for homogeneous $2$-intervals, which is tight,
and $\tau\leq (d^2-d+1)\nu$ for homogeneous $d$-intervals,
$d\geq 3$ \cite{Kaiser}. The reason for the worse bounds is that the
escape hypergraph needs no longer be $d$-partite, and so
F\"uredi's theorem \cite{Furedi-nunustar}
relating $\nu$ to $\nu^*$ gives a little worse
bound (for $d=2$, one uses a theorem of Lov\'asz instead,
asserting that $\nu^*\leq \frac 32\nu$ for any graph).

Sgall's construction \cite{Sgall} answered a problem raised by
Wigderson in 1985. The title of Sgall's paper
refers to a different, but essentially equivalent, formulation of
the problem dealing with labeled tournaments.

Alon \cite{Alon-trees} proved by the method
of \cite{Alon-piercing} that if $T$ is a tree and $\FF$ is
a family subgraphs of~$T$
with at most $d$ connected components, then
$\tau(\FF)\leq 2d^2\nu(\FF)$. More generally, he established 
a similar bound for the situation where $T$ is a graph of bounded
tree-width (on the other hand, if the tree-width of~$T$ is sufficiently
large, then one can find a system of connected subgraps of~$T$
with $\nu=1$ and $\tau$ arbitrarily large, and so the tree-width
condition is also necessary in this sense). A somewhat weaker bound for trees has been
obtained independently by Kaiser \cite{Kaiser-thesis}.

Strong results for piercing of $d$-trees, improving on Alon's results, 
were obtained by Berger \cite{Berger-d-trees},
based on a topological approach via KKM-type lemmas.
(For these see the references given above.)
\end{bibpar}

\begin{exs}
\item\label{exerc:intervals} We have claimed that for any family $\FF$ of intervals, 
it is well-known and easy to prove that $\nu(\FF)=\tau(\FF)$. 
Prove this!  
  
\item\label{ex:polytprod}
Let $P$ and $Q$ be convex polytopes. Show that there is
a bijection between the nonempty faces of
the Cartesian product $P\times Q$ and  all  the products
$F\times G$, where $F$ is a nonempty
face of~$P$ and $G$ is a nonempty face
of~$Q$.

\item
Show that the following ``Brouwer-like'' claim
resembling Lemma~\ref{l:polyt-surj}
is \emph{not} true: If $f\:B^n\to B^n$ is a continuous map
of the $n$-ball
such that the boundary of~$B^n$ is mapped surjectively onto
itself, then $f$ is surjective.

\item\label{ex:homog}
Prove the bound $\tau(\FF)\leq d^2\nu(\FF)$ for any
family of \emph{homogeneous} $d$-intervals (unions of~$d$
intervals on a single line). Hint: Follow the proof for
$d$-intervals above, but encode a candidate transversal $T$
by a point of a simplex (rather than  a product of simplices).
\end{exs}

\bigskip
\section{Evasiveness}
\label{sec:Evasiveness}

\subsection{A general model}

The idea of evasiveness comes from  the theory of complexity
of algorithms.
Evasiveness appears in different versions for graphs, digraphs
and bipartite graphs. We start with a general model
that contains them all.

\begin{Definition}
[Argument complexity of a set system; evasiveness]\label{d:ask}
In the following, we are concerned with a fixed and known set system $\FF\sse 2^E$,
and with the complexity of deciding whether some unknown set $A\sse E$
is in the set system. Here our ``model of computation'' is such that
\begin{description}
\itemsep=-3.5pt
\item[\quad given, and known,] is a set system $\FF\sse 2^E$, where 
$E$ is fixed, $|E|=m$.
\item[\quad]On the other hand, there is a
\item[\quad fixed, but unknown] subset $A \sse E$. 
\item[\quad]We have to 
\item[\quad decide] whether $A\in \FF$, using only 
\item[\quad questions] of the type ``Is $e\in A$?''
\end{description}
(It is assumed that we always get correct answers YES or NO.
We only count the \emph{number} of questions that are needed in order
to reach the correct conclusion: It is
assumed that it is not difficult to decide whether $e\in A$.
You can assume that some ``oracle'' that knows both $A$ and $\FF$ 
is answering.)

The \emph{argument complexity} $c(\FF)$ of the set system $\FF$ is 
the number of elements of the ground set $E$ that 
we have to test in the worst case---with the optimal strategy.

Clearly $0\le c(\FF) \le m$.
The set system $\FF$ is \emph{trivial} if $c(\FF)=0$:
then no questions need to be asked; this can only be the case if $\FF=\{\}$
or if $\FF=2^E$. Otherwise $\FF$ is \emph{non-trivial}.

The set system $\FF$ is \emph{evasive} if $c(\FF)=m$, that is, if even with an
optimal strategy one has to test all the elements of~$E$ in the worst case.
\end{Definition}
 
For example, if $\FF=\{\es\}$, then $c(\FF)=m$: If we again and again
get the answer
NO, then we have to test all the elements to be sure that $A=\es$.
So $\FF=\{\es\}$ is an evasive set system:
``being empty'' is an evasive set property.

\subsection{Complexity of graph properties}

\begin{Definition} [Graph properties]\label{d:gp}
Here we consider graphs on a fixed vertex set $V=[n]$.
Loops and multiple edges are excluded. Thus any graph 
$G=(V,A)$ is determined by its edge set~$A$, which is 
a subset of the set $E=\binom{n}{2}$ of ``potential edges.''

We identify a \emph{property} $\PP$ of graphs with the family of 
graphs that have the property~$\PP$, and thus with the set family
$\FF(\PP)\sse 2^E$ given by
\[
\FF(\PP)\assg\{A\sse E\sep \ ([n],A)\mbox{ has property }\PP\}.
\]
Furthermore, we will consider only graph properties that
are isomorphism invariant; that is, properties of abstract graphs
that are preserved under renumbering the vertices. 

A graph property is \emph{evasive} if the associated set system
is evasive, and otherwise it is \emph{non-evasive}.
\end{Definition}

With the symmetry condition of Definition~\ref{d:gp}, we would accept
``being connected'', ``being planar,'' 
``having no isolated vertices,'' and ``having even vertex degrees''
as graph properties.
However, ``vertex~$1$ is not isolated,'' ``$123$ is a triangle,'' and 
``there are no edges between odd-numbered vertices'' are not
graph properties.   

\begin{Examples}[Graph properties]
For the following properties of graphs on $n$ vertices we can
easily determine the argument complexity.
\begin{description}
\item[Having no edge:]
Clearly we have to check every single $e\in E$ in order to be sure
that it is not contained in $A$, so this property is evasive:
Its argument complexity is $c(\FF)=m=\binom{n}{2}$.
\item[Having at most \emph{k} edges:]
Let us assume that we ask questions, and the answer we get is YES 
for the first $k$ questions, and then we get NO answers
for all further questions, except for possibly the last one.
Assuming that $k<m$, this implies that the property is 
evasive. Otherwise, for $k\ge m$, the property is trivial.
\item[Being connected:]
This property is evasive for $n\ge2$. Convince yourself that
for any strategy, a sequence of ``bad'' answers can force you
to ask all the questions.
\item[Being planar:]
This property is trivial for $n\le 4$ but evasive for $n\ge5$.
In fact, for $n=5$ one has to ask all the questions
(in arbitrary order), and the answer will be $A\in\FF$ unless we get
a YES answer for all the questions---including the last one.
This is, however, not at all obvious for $n>5$: It was
claimed by Hopcroft \& Tarjan \cite{HopcroftTarjan},  
and proved by Best, Van Emde Boas \& Lenstra \cite[Example~2]{BEBL}
\cite[p. 408]{Bollobas}.  

\item[A large star:]
Let $\PP$ be the property of being a disjoint union of a
star $\Delta_{1,n-4}$ and an arbitrary graph on $3$ vertices, and let
$\FF$ be the corresponding set system.
\[
\input EPS/scorpion1.pstex_t
\]
Then $c(\FF)<\binom{n}{2}$ for $n\ge7$. For $n\ge12$ we can easily
see this, as follows. 
Test all the $\lfloor{n\over2}\rfloor\lceil{n\over2}\rceil$
edges $\{i,j\}$ with $i\le\lfloor{n\over2}\rfloor<j$. 
That way we will find exactly one vertex $k$ with at least 
$\lfloor{n\over2}\rfloor-3\ge 3$ neighbors
(otherwise property $\PP$ cannot be satisfied): That vertex $k$ has to be
the center of the star.
We test all other edges adjacent to~$k$: We must find that $k$ has
exactly $n-4$ neighbors. Thus we have identified three vertices
that are not neighbors of~$k$: At least one of the edges
between those three has not been tested. We test all other edges to
check that $([n],A)$ has property~$\PP$.
(This property was found by L. Carter \cite[Example~16]{BEBL}.)

\item[Being a scorpion graph:]
A \emph{scorpion graph} is an $n$-vertex graph that
has one vertex of degree~$1$ adjacent to a vertex of degree $2$
whose other neighbor has degree $n-2$.
We leave it as an (instructive!) exercise to check that 
``being a scorpion graph'' is not evasive if $n$ is large: In fact,
Best, van Emde Boas \& Lenstra \cite[Example~18]{BEBL} \cite[p.~410]{Bollobas}
have shown that $c(\FF)\le6n$.
\[
\input EPS/scorpion.pstex_t
\]

\end{description}
\end{Examples}

{From} these examples it may seem that 
most ``interesting'' graph properties are evasive.
In fact, many more examples of evasive graph properties can be found
in Bollob\'as \cite[Sect. VIII.1]{Bollobas}, alongside with
techniques to establish that graph properties are 
evasive, such as Milner \& Welsh's ``simple strategy'' 
\cite[p. 406]{Bollobas}.

Why is this model of interest?
Finite graphs (similarly for digraphs and bipartite graphs) 
can be represented in different types
of \emph{data structures} that are not at all equivalent for
algorithmic applications.
For example, 
if a finite graph is given by an \emph{adjacency list}, 
which for for every vertex lists the neighbors in some order, then one can
decide fast (``in linear time'') whether the graph is planar,
e.g.\ using an old algorithm of Hopcroft \& Tarjan \cite{HopcroftTarjan};
see also Mehlhorn \cite[Sect. IV.10]{Mehlhorn} and \cite{MehlhornMutzel}.
Note that such a planar graph has at most $3n-6$ edges (for $n\ge3$).
  
However, assume that a graph is given in terms of its adjacency matrix
\[
M(G)\ \ =\ \ \big(m_{ij}\big)_{1\le i,j\le n}\ \ \in\ \{0,1\}^{n\times n},
\]
where $m_{ij}=1$ means that $\{i,j\}$ is an edge of~$G$,
and $m_{ij}=0$ says that $\{i,j\}$ is not an edge.
Here $G$ is faithfully represented by the set of all $\binom{n}{2}$ superdiagonal
entries (with $i<j$).
Then one possibly has to inspect a large part of the
matrix until one {has enough information
to decide whether the graph in question is planar}. 
In fact, if $\FF\sse 2^E$ is the set system corresponding to all planar graphs,
then $c(\FF)$ is exactly the number of superdiagonal matrix entries
that every algorithm for planarity testing has to
inspect in the worst case.

The statement that ``being planar'' is evasive (for $n\ge5$) thus
translates into the fact that every planarity testing algorithm
that starts from an adjacency matrix needs to read at least
$\binom{n}{2}$ bits of the input, and hence its running time
is bounded from below by $\binom{n}{2}=\Omega(n^2)$.
This means that such an algorithm---such as the one considered
by Fisher \cite{Fisher-planar}---cannot run in linear time,
and thus cannot be efficient.

\begin{Definition} 
[Digraph properties; bipartite graph properties]~
\begin{compactenum}[(1)]
 \item
For digraph properties we again use the fixed vertex set $V=[n]$.
Loops and parallel edges are excluded, but
anti-parallel edges are allowed. Thus any digraph 
$G=(V,A)$ is determined by its arc set~$A$, which is 
a subset of the set $E'$ of all $m\assg n^2-n$ ``potential arcs''
(corresponding to the off-diagonal entries of an $n\times n$ 
adjacency matrix).

A \emph{digraph property} is a property of digraphs
$([n],A)$ that is invariant under relabelling of the vertex set.
Equivalently, a digraph property is a family of arc sets $\FF\sse 2^{E'}$
that is symmetric under the action of~$\Symm_n$ that acts by
renumbering the vertices (and renumbering all arcs correspondingly).
A digraph property is \emph{evasive} if the associated set system
is evasive, otherwise it is \emph{non-evasive}.
 \item
For bipartite graph properties we use a fixed vertex set 
$V\uplus W$ of size $m+n$, and use $E''\assg V\times W$ 
as the set of potential edges.
A \emph{bipartite graph property} is a property of graphs
$(V\cup W, A)$ with $A\sse E''$ that is preserved under 
renumbering the vertices in $V$, and also under permuting the
vertices in~$W$. Equivalently, a bipartite graph property on $V\times W$ is
a set system $\FF\sse 2^{V\times W}$ that is stable under the
action of the automorphism group
$\Symm_n\times \Symm_m$ that acts transitively on~$V\times W$.
\end{compactenum}
\end{Definition}

\begin{Examples}[Digraph properties]
For the following digraph properties on $n$ vertices we can
determine the argument complexity.
\begin{description}
\item[Having at most \emph{k} arcs:]
Again, this is clearly evasive with $c(\FF)=m$ if $k<m=n^2-n$,
and trivial otherwise.
\item[Having a sink:]
A \emph{sink} in a digraph on $n$ vertices is a 
vertex $k$ for which all arcs going into $k$ are present,
but no arc leaves $k$, that is, a vertex of
out-degree $\delta^+(v)=0$, and in-degree $\delta^-(v)=n-1$.
Let $\FF$ be the set system of all 
digraphs on $n$ vertices that have a sink.
It is easy to see that $c(\FF)\le 3n-4$. In particular, for $n\ge3$
``having a sink'' is a non-trivial but non-evasive digraph property.

In fact, if we test whether $(i,j)\in A$, then either
we get the answer YES, then $i$ is not a sink, or  
we get the answer  NO, then $j$ is not a sink.
So, by testing arcs between pairs of vertices that ``could be sinks,''
after $n-1$ questions we are down to one single ``candidate sink''~$k$. 
At this point at least one arc adjacent to~$k$ has been tested.
So we need at most $2n-3$ further questions to test whether it is a sink.
\end{description}
\end{Examples}
 
In the early 1970's Arnold L. Rosenberg had conjectured that all 
non-trivial digraph properties have qua\-dra\-tic argument complexity,
that is, that there is a constant $\gamma>0$ such that for all
non-trivial properties of digraphs on $n$ vertices
one has $c(\FF) \ge \gamma n^2$. However, St{\aa}l Aanderaa found the
counter-example (for digraphs) of ``having a sink'' \cite[Example~15]{BEBL} \cite[p. 372]{RV1}. We have also seen that ``being a scorpion graph''
is a counter-example for graphs.

Hence Rosenberg modified the conjecture: 
At least all \emph{monotone} graph properties, that is,
properties that are preserved under deletion of edges, should have
quadratic argument complexity. This is the statement of the 
\emph{Aanderaa--Rosenberg conjecture} \cite{Rosenberg}. 
Richard Karp considerably sharpened the statement, as follows. 

\begin{Conjecture} [The evasiveness conjecture]
Every non-trivial monotone graph property or digraph property is evasive.
\end{Conjecture}

We will prove this below for graphs and digraphs in the
special case when $n$ is a prime power;
from this one can derive
the Aanderaa--Rosenberg conjecture, with $\gamma\approx{1\over4}$.
Similarly, we will prove that monotone properties of bipartite graphs on a
fixed ground set $V\cup W$  are evasive (without any restriction
on $|V|=m$ and $|W|=n$). 
However, we first return to the more general setting of set systems.

\subsection{Decision trees}

Any strategy to determine whether an (unknown) set  $A$ is contained
in a (known) set system~$\FF$---as in 
Definition~\ref{d:ask}---can be represented in terms 
of a decision tree of the following form.
 
\begin{Definition}\rm
A \emph{decision tree} is a rooted, planar, binary tree
whose leaves are labelled ``YES'' or ``NO,'' and whose internal
nodes are labelled by questions (here they are of the type ``$e\in A?$'').
Its edges are labelled by answers: We will represent them so
that the edges labelled ``YES'' point to the right child, and the
``NO'' edges pointing to the left child.

A \emph{decision tree for $\FF\sse 2^E$} is a decision tree such 
that starting at the root with an arbitrary $A\sse E$,
and going to the right resp.\ left child depending on whether
the question at an internal node we reach has answer YES or NO,
we always reach a leaf that correctly answers the question
``$A\in\FF?$''.
\[
\input{EPS/inner2.pstex_t}
\] 
The root of a decision tree is at \emph{level}~$0$, and the children
of a node at level $i$ have level $i+1$.
The \emph{depth} of a tree is the greatest $k$ such that the tree
has a vertex at level~$k$ (a leaf).

We assume (without loss of generality) that the trees 
we consider correspond to strategies where we never ask the same question twice.

A decision tree for $\FF$ is \emph{optimal} if it has the smallest depth 
among all decision trees for~$\FF$, that is, if it 
leads us to ask the smallest number of questions
for the worst possible input.
\end{Definition}

Let us consider an explicit example. 
\[
\input{EPS/k3d.pstex_t}
\]
The following figure represents an optimal algorithm
for the ``sink'' problem on digraphs with $n=3$ vertices.
This has a ground set 
$E=\{12,21,13,31,23,32\}$ of size $m=6$.

\[
\hspace*{-10pt}\input{EPS/tree2.pstex_t}
\]
The algorithm first asks, in the root node at {level~$0$},
whether $12\in A$. 
In case the answer is YES (so we know that $1$ is not a sink),
it branches to the right,
leading to a question node at level~$1$ that asks whether $23\in A?$,
etc.
In case the answer to the question $12\in A?$ is NO  
(so we know that $2$ is not a sink),
it branches to the left,
leading to a question node at level~$1$ that asks whether $13\in A?$,
etc.  

For every possible input $A$ (there are $2^6=32$ different ones),
after two questions we have identified a unique ``candidate sink'';
after not more than $5$ question nodes one arrives at a leaf node
that correctly answers the question whether the graph $(V,A)$ has a sink
node: YES or NO.
(The number of the unique candidate is noted next to each node at 
level~$2$.)

For each node (leaf or inner) of level $k$, there are
exactly $2^{m-k}$ different inputs that lead to this node. 
This proves the following lemma.

\begin{Lemma} \label{l:evas-char}
The following are equivalent:
\begin{compactitem}[ $\bullet$]
\item $\FF$ is non-evasive. 
\item The optimal decision trees $T_\FF$ for $\FF$ have depth smaller than~$m$.
\item Every leaf of an optimal decision tree $T_\FF$ is reached by
at least two distinct inputs.
\end{compactitem}
\end{Lemma}
 
\begin{Corollary} 
If $\FF$ is non-evasive, then $|\FF|$ is even.
\end{Corollary} 

This can be used to show, for example, that
the directed graph property ``has a directed cycle'' is evasive
\cite[Example~4]{BEBL}.

Another way to view a (binary) decision tree algorithm is as follows.
In the beginning, we do not know anything about the set $A$, so we can view
the collection of possible sets as the complete boolean algebra of
all $2^m$ subsets of~$E$.%

In the first node (at ``level~$0$'')
we ask a question of the type ``$e\in A?$''; this induces a subdivision of
the boolean algebra into two halves, depending on whether we get
answer YES or NO. If you think of the boolean algebra as a partially ordered set
(indeed, a lattice), then each of the halves is an interval of length $m-1$
of the boolean algebra $(2^E,\sse)$.
If you prefer to think of it as a rendition of the $m$-dimensional hypercube,
then the halves are subcubes of codimension $1$, containing all the vertices
of two opposite facets.  
  
At level~$1$ we ask a new question, depending
on the outcome of the first question. Thus we 
\emph{independently} bisect the two halves of level $0$,
getting four pieces of the boolean algebra, all of the same size.
\[
\input EPS/boolean5.pstex_t  
\]  
This process is iterated. It stops---as we do not need to ask
a further question---on parts that we create that
either contain only sets that are in $\FF$ (this yields a YES-leaf) 
or that contain only sets not in $\FF$ (corresponding to NO-leaves).

Thus the final result is a special type of partition of the
boolean algebra into intervals. 
Some of them are YES intervals, containing only sets of~$\FF$, 
all the others are NO-intervals, containing no sets from~$\FF$.
If the property in question is monotone, then the 
union of the YES intervals (i.\,e., the set system~$\FF$) forms an \emph{ideal}
in the boolean algebra, that is, a ``down-closed'' set such that with any 
set that it contains it must also contain all its subsets.

Let $p_\FF(t)$ be the generating function for the 
set system $\FF$, that is, the polynomial 
\[
p_\FF(t)\assg\sum_{A\in\FF} t^{|A|} \ \ =\ \  
f_{-1}+tf_0+t^2f_1+t^3f_2+\dots.
\]
where $f_i= |\{A\in\FF\sep |A|=i+1\}|$.

\begin{Proposition} 
\[
(1+t)^{m-c(\FF)}\ \ \big|\ \ p_{\FF}(t).
\]
\end{Proposition}

\proof
Consider one interval $\II$ in the partition of~$2^E$
that is induced by any optimal algorithm for $\FF$.
If the  leaf, at level $k$, corresponding to the interval is reached through
a sequence of~$k_Y$ YES-answers and $k_N$ NO-answers (with $k_Y+k_N=k$),
then this means that there are sets $A_Y\sse E$ with $|A_Y|=k_Y$
and  $A_N\sse E$ with $|A_N|=k_N$, such that 
\[
\II\ \ =\ \ \{A\sse E\sep A_Y\sse A\sse E\sm A_N\}.
\]
In other words, the interval $\II$ contains all sets that give YES-answers
when asked about any of the $k_Y$~elements of~$A_Y$, NO-answers
when asked about any of the $k_N$~elements of~$A_N$, while
the $m-k_Y-k_N$ elements of $E\sm(A_Y\cup A_N)$ may or may not be contained
in~$A$.
Thus the interval $\II$ has size $2^{m-k_Y-k_N}$, and its counting
polynomial is 
\[
p_\II(t)\assg\sum_{A\in\II} t^{|A|} \ \ =\ \  t^{k_Y}(1+t)^{m-k_Y-k_N}.
\]
Now the complete set system~$\FF$ is a disjoint union of the intervals~$\II$,
and we get
\[
p_\FF(t)\ \ =\ \ \sum_{\II} p_\II(t).
\]
In particular, for an optimal decision tree we have
$k_Y+k_N=k\le c(\FF)$ and thus $m-c(\FF)\le m-k_Y-k_N$
at every leaf of level~$k$, which means that
all the summands $p_{\II}(t)$ have a common factor of $(1+t)^{m-c(\FF)}$.
\endproof

\begin{Corollary} \label{c:euler-char} 
If $\FF$ is non-evasive, then $|\FF^{even}| = |\FF^{odd}|$,
that is, 
\[-f_{-1}+f_0-f_1+f_2\mp\cdots = 0.\]
\end{Corollary} 

\proof
Use Lemma~\ref{l:evas-char}, and put $t=-1$.
\endproof

We can now draw the conclusion, based only on simple counting, that most set families are evasive. This cannot of course be used to settle any specific cases,
but it can at least make the various evasiveness conjectures 
seem more plausible.

\begin{Corollary}  
Asymptotically, almost all set families $\FF$ are evasive.
\end{Corollary} 

\proof
The number of set families $\FF\subseteq 2^E$ such that
\[
\# \{A\in \FF \mid \# A \mbox{ odd} \}  = 
\# \{A\in \FF \mid \# A \mbox{ even} \}  =k
\]
is ${\binom{2^{m-1}}{k}}^2$. Hence, using Stirling's estimate of factorials,
\[ \mbox{ Prob ($\FF$ non-evasive)} \, \le\, 
\frac{\sum_{k=0}^{2^{m-1}} {\binom{2^{m-1}}{k}}^2}{2^{2^m}}\, = \,
\frac{\binom{2^{m}}{2^{m-1}}}{2^{2^m}} \, \sim\, \frac{1}{\sqrt{\pi2^{m-1}}}
\rightarrow 0,
\]
as $m\rightarrow \infty$.
\endproof

\begin{Conjecture}[The ``Generalized Aanderaa--Rosenberg Conjecture'', Rivest \& Vuillemin \cite{RV2}] \label{c:GAR}
If $\FF\sse 2^E$, with symmetry group $G\sse \Symm_E$ that 
is transitive on the ground set $E$, and if 
$\es\in\FF$ but $E\notin\FF$, then $\FF$ is evasive.
\end{Conjecture}

Note that for this it is \emph{not} assumed that
$\FF$ is monotone. However, the assumption that $\es\in\FF$ but $E\notin\FF$
is satisfied neither by ``being a scorpion'' nor by ``having a sink.''

\begin{Proposition}[\rm Rivest \& Vuillemin \cite{RV2}]\label{p:RV}
The Generalized Aanderaa--Rosenberg 
Conjecture~\ref{c:GAR}
holds if the size of the ground set is a prime power, $|E|=p^t$.
\end{Proposition}

\proof
Let $\OO$ be any $k$-orbit of~$G$, that is, a collection of
$k$-sets $\OO\sse\FF$ on which $G$ acts transitively.
While every set in $\OO$ contains $k$ elements $e\in E$, we know from 
transitivity that every element of~$E$ is contained in the same
number, say $d$, of sets of the orbit~$\OO$.
Thus, double-counting the edges of the bipartite graph on the vertex
set $E\uplus\OO$ defined by ``$e\in A$'' (displayed in the figure below)
we find that $k|\OO|= d|E| = d p^t$.
Thus for $0<k<p^t$ we have that $p$ divides $|\OO|$,
while $\{\es\}$ is one single ``trivial'' orbit of size~$1$,
and $k=p^t$ doesn't appear.
Hence we have 
\[
-f_{-1}+f_0-f_1+f_2\mp\cdots \equiv -1 \bmod p,
\]
which implies evasiveness by Corollary~\ref{c:euler-char}.
\endproof

\[
\input EPS/orbit.pstex_t
\]

\begin{Proposition}[Illies \cite{Illies}]\label{e:illies}
The Generalized Aanderaa--Rosenberg Conjecture~\ref{c:GAR}
fails for $n=12$.
\end{Proposition}

\proof
Here is Illies' counterexample: Take 
$E=\{1,2,3,\dots,12\}$, and let the cyclic group 
$G=\Z_{12}$ permute the elements of~$E$ with the obvious cyclic action.

Take $\FF_I \sse 2^E$ to be the following system of sets
\begin{itemize}\itemsep=0pt
\item $\es$, so we have $f_{-1}=1$
\item $\{1\}$ and all images under $\Z_{12}$,
that is, all singleton sets: $f_0=12$,
\item $\{1,4\}$ and $\{1,5\}$ and all images under $\Z_{12}$, 
so $f_1=12+12=24$,
\item $\{1,4,7\}$ and $\{1,5,9\}$ and all their $\Z_{12}$-images, 
so $f_2=12+4=16$,
\item $\{1,4,7,10\}$ and their $\Z_{12}$-images, so $f_3=3$.
\end{itemize}

An explicit decision tree of depth~$11$ for this $\FF_I$
is given in our figure below. 
Here the \emph{pseudo-leaf} ``YES(7,10)'' denotes
a decision tree where we check all elements $e\in E$ that 
have not been checked before, other than the elements $7$ and $10$.
If none of them is contained in $A$, then the answer is YES
(irrespective of whether $7\in A$ or $10\in A$), 
otherwise the answer is~NO. The fact that two elements need not be
checked means that this branch of the decision tree denoted by this
``pseudo-leaf'' does not go beyond depth~$10$.
Similarly, a pseudo-leaf of the type ``YES(7)'' represents
a subtree of depth~$11$.

Thus the following figure completes the proof. Here dots denote subtrees that are analogous to the ones just above.
\endproof

\[\hspace{-20pt}
\input EPS/illies3.pstex_t
\]
Note, however, that Illies' example is not monotone: For example, we have
$\{1,4,7\}\in\FF_I$, whereas $\{1,7\}\notin\FF_I$.

\subsection{Monotone systems}

We now concentrate on the case where $\FF$ is closed under taking subsets,
that is, $\FF$ is an abstract simplicial complex, which we also denote 
by $\Delta\assg\FF$. In this setting, the symmetry group
acts on~$\Delta$ as a group of simplicial homeomorphisms.
If $\FF$ is a graph or digraph property, then this means
that the action of~$G$ is transitive on
the vertex set $E$ of~$\Delta$, which corresponds to the
edge set of the graph in question. Again we denote the cardinality of the
ground set (the vertex set of~$\Delta$) by $|E|=m$.

A complex $\Delta\sse 2^E$ is a \emph{cone} if it has 
a vertex $v$ such that $A\cup\{v\}$ is a face of~$\Delta$ for any face $A\in\Delta$.
For example, every simplex $\Delta= 2^E$ is a cone, but also every
star graph $K_{m,1}$, considered as a simplicial complex of dimension $1$, is a cone.

A complex $\Delta\sse 2^E$ is \emph{collapsible} if it can be reduced
to a one-point complex (equivalently, to a simplex) by steps of the form 
\[
\Delta\to \Delta\sm \{A\in \Delta\sep A_0\sse A\sse A_1\}
\]
$\es\subset A_0\subset A_1$ are faces of~$\Delta$ with $\es\neq A_0\neq A_1$, 
where $A_1$ is the \emph{unique} maximal element of~$\Delta$ that contains $A_0$.
For example, every tree, considered as a simplicial complex of dimension~$1$, 
is collapsible.

Our figure illustrates a sequence of collapses that reduce
a $2$-dimensional complex to a point. In each case the face $A_0$ 
that is contained in a unique maximal face is drawn fattened.
\[
\includegraphics[width=12cm]{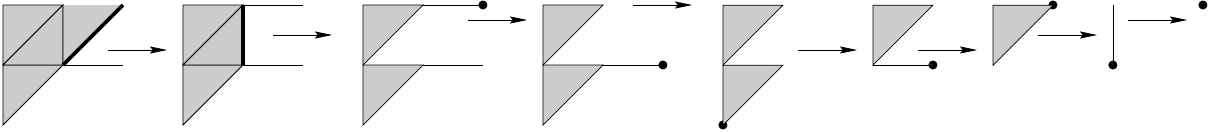}
\]

\begin{Theorem}\label{t:collapsible}
We have the following implications:
\\ 
$\Delta$ is a cone $\ \Lra\ $
$\Delta$ is non-evasive $\ \Lra\ $
$\Delta$ is collapsible $\ \Lra\ $
$\Delta$ is contractible.
\end{Theorem}

\proof
The first implication is clear: For a cone we don't have to test the apex~$e_0$
in order to see whether a set $A$ is a face of~$\Delta$, since
$A\in \Delta$ if and only if $A\cup\{e_0\}\in \Delta$.
The third implication is easy topology: One can write down 
explicit deformation retractions.
The middle implication we will derive from the following claim,
which uses the notion of a \emph{link} of a vertex $e$ in a simplicial
complex $\Delta$: This is the complex $\Delta/e$ formed by all faces $A\in\Delta$
such that $e\notin A$ but $A\cup\{e\}\in\Delta$.

\textbf{Claim.}
\emph{$\Delta$ is non-evasive if and only if either
$\Delta$ is a simplex, or it is not a simplex but it has a vertex $e$ such that
both the deletion $\Delta\sm e$ and the link $\Delta/e$ are non-evasive.}

Let us first verify this claim:
If no questions need to be asked (that is, if $c(\Delta)=0$), then $\Delta$ is a simplex.
Otherwise we have some $e$ that corresponds to
the first question to be asked by an optimal algorithm. 
If one gets a YES answer, then the problem
is reduced to the link~$\Delta/e$, since the faces $B\in \Delta/e$ 
correspond to the faces $A=B\cup\{e\}$ of~$\Delta$ for which $e\in A$.
In the case of a NO-answer the problem similarly reduces to 
the deletion~$\Delta\sm e$.

Now let us return to the proof of Theorem~\ref{t:collapsible},
where we still have to verify that
``$\Delta$ is non-evasive $\Lra$ $\Delta$ is collapsible.''
We use induction on the number of faces of~$\Delta$.

If $\Delta$ is not a simplex, then by the Claim
it has a vertex $e$ such that the link~$\Delta/e$ and the deletion $\Delta\sm e$
are collapsible.
If the link is a simplex, then deletion of~$e$ is a collapsing step
$\Delta\rightarrow\Delta\sm e$, where $\Delta\sm e$ is collapsible,
so we are done by induction.

If the link is not a simplex, then it has faces $\es\subset A_0\subset A_1$ such that $A_1$ is the unique
maximal face in the link that contains $A_0$.
This means that $\Delta$ has faces $\{e\}\subset A_0\cup\{e\}\subset A_1\cup\{e\}$ such that 
$A_1\cup\{e\}$ is the unique
maximal face in~$\Delta$ that contains $A_0\cup\{e\}$. 
In this way any collapsing step in the link $\Delta/e$
yields a collapsing step in $\Delta$, and again we are done by induction.
\endproof

\subsection{A topological approach}\label{sec4}

The following simple lemma provides the step from the topological
fixed point theorems for complexes to combinatorial information.

\begin{Lemma}\label{l:transitive-fixed}
If a (finite) group $G$ acts vertex-transitively on a finite complex~$\Delta$
with a fixed point, then $\Delta$ is a simplex.
\end{Lemma}

\proof
If $V\assg\{v_1,\dots,v_n\}$ is the vertex set of~$\Delta$, then
any point $x\in \polyh{\Delta}$ has a unique representation of the form
\[
x\ \ =\ \ \sum_{i=1}^n \la_i\,v_i ,
\]
with $\la_i\ge0$ and $\sum_{i=1}^n\la_i=1$.
If the group action, with
\[
gx\ \ =\ \ \sum_{i=1}^n \la_i\,gv_i,
\]
is transitive, then this means that for every $i,j$ there is some
$g\in G$ with $gv_i=v_j$. 
Furthermore, if $x$ is a fixed point, then we have $gx=x$ for all $g\in G$,
and hence we get $\la_i=\la_j$ for all $i,j$. From this we derive
$\la_i={1\over n}$ for all~$i$. Hence we get
\[
x\ \ =\ \ {1\over n}\sum_{i=1}^n \,v_i 
\]
and this is a point in~$\polyh{\Delta}$ only if $\Delta$ is the 
complete simplex with vertex set~$V$.
\smallskip

Alternatively: The fixed point set of any group action
is a subcomplex of the
barycentric subdivision, by Lemma~\ref{l:fps}.
Thus a vertex $x$ of the fixed point complex
is the barycenter of a face $A$ of~$\Delta$. Since $x$ is fixed by the
whole group, so is its support, the set $A$.
Thus vertex transitivity implies that $A=E$, and $\Delta=2^E$.
\endproof

\begin{Theorem} [The Evasiveness Conjecture for prime powers:
Kahn, Saks \& Sturtevant \cite{KSS}] 
All monontone non-trivial graph properties and digraph properties for 
graphs on a prime power number of vertices $|V|=q=p^t$ are evasive.
\end{Theorem}

\proof
We identify the fixed vertex set $V$ with $\GF(q)$.
Corresponding to a non-evasive monotone non-trivial graph property we
have a non-evasive complex $\Delta$ on a set $E=\binom{V}{2}$ 
of~$\binom{q}{2}$ vertices.
By Theorem~\ref{t:collapsible} $\Delta$ is collapsible and hence \emph{$\Z_p$-acyclic},
that is, all its reduced homology groups with $\Z_p$-coefficients vanish.

The symmetry group of~$\Delta$ includes the symmetric group~$\Symm_q$,
but we take only the subgroup of all ``affine maps''  
\begin{eqnarray*}
G  & \assg   &  \{x\mapsto ax+b \sep a,b\in \GF(q),\ a\neq0\},\\
\noalign{\noindent and its subgroup}
P  &  \assg &  \{x\mapsto  x+b \sep b\in \GF(q)\} 
\end{eqnarray*}
that permute the vertex set~$V$, and (since we are considering graph properties)
extend to an action on  
the vertex set $E=\binom{V}{2}$ of~$\Delta$.  
Then we can easily verify the following facts:
\begin{compactitem}[ $\bullet$] 
\item $G$ is doubly transitive on $V$, and hence induces a vertex transitive
group of symmetries of the complex $\Delta$ on the vertex set~$E=\binom{V}{2}$ 
(interpret $\GF(q)$ as
a $1$-dimensional vector space, then any (ordered) pair of distinct
points can be mapped to any other such pair by an affine map on the line);
\item $P$ is a $p$-group (of order $p^t=q$);
\item $P$ is the kernel of the homomorphism that maps
$(x\mapsto ax+b)$ to $a\in \GF(q)^*$, the multiplicative group of~$\GF(q)$,
and thus a normal subgroup of~$G$;
\item $G/P\cong \GF(q)^*$ is cyclic (this is known from your algebra class).
\end{compactitem} 
Taking these facts together, we have verified all the 
requirements of Oliver's fixed point theorem, as provided in the Appendix as Theorem~\ref{t:oliver}. 
Hence
$G$ has a fixed point on $\Delta$, and by Lemma~\ref{l:transitive-fixed}
$\Delta$ is a simplex, and hence the corresponding (di)graph property is
trivial.
\endproof

From this one can also deduce---with a lemma due to
Kleitman \& Kwiatowski \cite[Thm. 2]{KleitmanKwiatowski}---that
every non-trivial monotone graph property on $n$ vertices
has complexity at least $n^2/4+o(n^2) = m/2 + o(m)$.
(For the proof see \cite[Thm. 6]{KSS}.)
This establishes the Aanderaa--Rosenberg Conjecture.  
On the other hand, the
Evasiveness Conjecture is still an open problem for every
$n\ge10$ that is not a
prime power. Kahn, Saks \& Sturtevant \cite[Sect. 4]{KSS} report that they
verified it for $n=6$.

The following treats the bipartite version of the 
Evasiveness Conjecture. Note that in the
case where $mn$ is a prime power it follows from Proposition~\ref{p:RV}.

\begin{Theorem} [The Evasiveness Conjecture for bipartite graphs,
Yao~\cite{Yao}]
All monotone non-trivial bipartite graph properties are evasive.
\end{Theorem}

\proof
The ground set now is $E=V\times W$, where
any monotone bipartite graph property is represented by 
a simplicial complex $\Delta\sse 2^E$.

An interesting aspect of Yao's proof is that it does not
use a vertex transitive group. In fact, let the cyclic group
$G\assg \Z_n$ act by cyclically permuting the vertices in $W$,
while leaving the vertices in $V$ fixed.
The group $G$ satisfies the assumptions of Oliver's Theorem~\ref{t:oliver},
with $P=\{0\}$. It acts on the complex $\Delta$ which is acyclic by
Theorem~\ref{t:collapsible}. Thus we get from Oliver's Theorem that
the fixed point set $\Delta^G$ is acyclic.
This fixed point set is not a subcomplex of~$\Delta$ (it does not contain
any vertices of~$\Delta$), but it is a subcomplex of the
order complex $\Delta(\Delta)$, which is the barycentric subdivision of~$\Delta$
(Lemma~\ref{l:fps}).

The bipartite graphs that are fixed under $G$ are those for which 
every vertex in $V$ is adjacent to none, or to all, of the vertices in~$W$;
thus they are complete bipartite graphs of the type $K_{k,n}$ for 
suitable $k$. Our figure illustrates this for the case where 
$m=6$, $n=5$, and $k=3$.
\[
\input EPS/Kkn.pstex_t
\]
Monotonicity now implies that the fixed graphs under $G$ are \emph{all} the
complete bipartite graphs of type $K_{k,n}$ with $0\le k\le r$ for some
$r$ with $0\le r< m$. (Here $r=m$ is impossible, since then $\Delta$ would be
a simplex, corresponding to a trivial bipartite graph property.)

Now we observe that $\Delta^G$ is the order complex (the barycentric
subdivision) of a different complex, namely of the complex whose
vertices are the complete bipartite subgraphs $K_{1,n}$, and whose
faces are \emph{all} sets of at most $r$ vertices.

Thus $\Delta^G$ is the barycentric subdivision of the $(r-1)$-dimensional skeleton
of an $(m-1)$-dimensional simplex. In particular, this space is not
acyclic. Even its reduced Euler characteristic, which can be
computed to be $(-1)^{r-1}\binom{m-1}{r}$, does not vanish.
\endproof

\noindent 
We have the following sequence of implications:   \[
\mbox{%
non-evasive\textsuperscript{(1)} $\Lra$ 
collapsible\textsuperscript{(2)} $\Lra$
contractible\textsuperscript{(3)} $\Lra$
$\Q$-acyclic\textsuperscript{(4)} $\Lra$ $\chi =1$\textsuperscript{(5)},}
\] 
which corresponds to a sequence of conjectures:
\\[2mm]
\textbf{Conjecture(\emph{k}):}
\emph{Every vertex-homogeneous simplicial complex with property~$(k)$ is a simplex.}
\\[2mm]
The above implications show that 
\begin{small}%
\[
\mbox{Conj.\,(5)}  \Lra 
\mbox{Conj.\,(4)}  \Lra 
\mbox{Conj.\,(3)}  \Lra 
\mbox{Conj.\,(2)}  \Lra 
\mbox{Conj.\,(1)}  \Lra
\begin{array}{l}
    \mbox{\small Evasiveness}\\ \mbox{\small Conjecture}    
\end{array}
\]
\end{small}%
Here Conjecture (5) is \emph{true} for a prime power number of vertices,
by Theorem~\ref{p:RV}.

However, Conjectures (5) and (4) fail for $n=6$: 
A counterexample is
provided by the six-vertex triangulation of the real projective
plane (see \cite[Section~5.8]{Mat-top}).
Even Conjectures (3) and possibly (2) fail for $n=60$:
a counterexample by Oliver (unpublished), of dimension~$11$,
is based on the group $A_5$; see Lutz~\cite{lutz02:_examp_z}.

So, it seems that Conjecture (1)---the monotone version of the
Generalized Aanderaa--Rosenberg 
Conjecture~\ref{c:GAR}---may be the right generality to prove, 
even though its
non-monotone version fails by Proposition~\ref{e:illies}. 

\subsection{Quillen's conjecture}
\label{sec4.5}

In this final section we briefly comment on a well-known conjecture of
Daniel Quillen from 1978 concerning finite groups. Upon
first sight it seems very remote from 
the topic of evasiveness that we have just discussed,
but under the surface one finds some surprising similarities. 

In this section we assume familiarity with basic finite group theory, 
and with the topology of order complexes.

A finite group is a \emph{$p$-group} if its order is a power of the 
prime number $p$. A subgroup of a finite group $G$ is a \emph{$p$-Sylow subgroup} if it is a maximal 
$p$-group. The number $n_p$  of $p$-Sylow subgroups of~$G$ is called
the \emph{$p$-Sylow number} of~$G$.
 
Let $G$ be a finite group and $p^e$ a prime power such that
 $|G| =p^e m$ and $p$ does not divide $m$. 
Here are some well known properties.
\begin{enumerate}
\item There exists a $p$-Sylow subgroup of~$G$ of order 
$p^{e}$.
\item  Any two  $p$-Sylow subgroups of~$G$ are conjugate to each other. 
\item  $n_p(G) \equiv 1 \mod{p}$.
\end{enumerate}
These statements are the familiar \emph{Sylow theorems}, 
the first substantial results in most treatises on group theory. 

For a finite group $G$ and a prime number $p$
dividing its order, let $L_p(G)$
denote the poset of non-trivial $p$-subgroups of~$G$, ordered by inclusion.
This is a ranked poset, the maximal elements 
of which are the $p$-Sylow subgroups. It
becomes a lattice if one adds new bottom and top elements.

In 1978 Quillen published the following conjecture \cite{Quillen1978}, which
in a surprising way connects a topological condition with an algebraic one.

\begin{Conjecture} [Quillen's conjecture]
$L_p(G)$ is contractible if and only if $G$ has a non-trivial normal $p$-subgroup.
\end{Conjecture}

Here $L_p(G)$ refers to the order complex, whose simplices are the 
totally ordered chains $x_0 < x_1 < \dots < x_d$ of~$L_p(G)$.
The ``if'' direction, which is very easy, was proved by Quillen,
and he proved the ``only if'' direction for the case of solvable groups.
The conjecture has since then been verified in many cases, but 
the general case is still wide open. 

In the previous section we considered an array of conjectures,
among them this one: 
\smallskip
 
\noindent\textbf{Conjecture (3):}
\emph{Every vertex-homogeneous contractible simplicial complex is a simplex.}
\smallskip

This conjecture turns out to be relevant both for evasiveness and for $p$-subgroups:
\begin{quote}
	 Conjecture (3) $\Longrightarrow$ Evasiveness Conjecture,
	 
	 Conjecture (3) $\Longrightarrow$ Quillen's Conjecture.
\end{quote} 
However, Conjecture (3) is false. It was mentioned in the 
previous section that
counterexamples on $60$ vertices are known.
So, why spend time on discussing it? We believe that it is
nevertheless instructive to see in which way
Conjecture (3) is relevant for Quillen's Conjecture.
It is conceivable that progress for one of the Evasiveness Conjecture and
the Quillen Conjecture can lead to progress for the other.
 
\begin{Proposition}
Conjecture $(3)$ $\Longrightarrow$ Quillen's Conjecture
\end{Proposition}

\begin{proof}
Suppose that  $L_p(G)$ is contractible. We are to prove that $G$ has a non-trivial
normal $p$-subgroup. 

Define the auxiliary \emph{Sylow complex} $\mathrm{Syl}_p(G)$ this way: The vertices are the $p$-Sylow subgroups
of~$G$. A collection of such subgroups form a simplex (or, face) of 
$\mathrm{Syl}_p(G)$
 if their intersection is nontrivial (not just the identity).
This is clearly a simplicial complex.

An application of the nerve theorem (or the crosscut theorem), see Björner \cite[p.~1850]{Bjorner-Handbook},
shows that these two complexes are of same homotopy type:
\[
 \mathrm{Syl}_p(G) \sim L_p(G)
\]
The group $G$ acts by conjugation on the vertex set of
$\mathrm{Syl}_p(G)$, and by the second Sylow theorem this action is  transitive.
So, $\mathrm{Syl}_p(G)$ is a vertex-homogeneous and contractible complex.
Conjecture (3) then implies that $\mathrm{Syl}_p(G)$ is a big simplex.
This means precisely that the intersection of \emph{all} $p$-Sylow subgroups is
non-trivial and is a fixed point under the action. Hence this is a non-trivial
normal $p$-subgroup.
\end{proof}

Following along the reasoning in this proof can help to verify
the Quillen conjecture in some special cases,
such as this.

\begin{Proposition}
If  $n_p= q^{e}$, that is, if the number of $p$-Sylow subgroups
 is the power of some prime number $q$,
then $G$ satisfies the Quillen conjecture.
\end{Proposition}

Here the Rivest--Vuillemin Theorem~\ref{p:RV} is relevant. In fact, with this and 
Conjecture (5) a sharper version of the Quillen conjecture can be obtained
in the case when $n_p= q^{e}$,
using trivial Euler characteristic instead of contractibility.
We leave further thoughts and experiments in this direction to the reader.

\begin{bibpar}
	The classical textbook account on evasiveness, from the Graph Theory 
point of view, is in Bollobas \cite[Chap. VIII]{Bollobas}. 

A textbook account from a Topological Combinatorics point-of-view was recently
given in de Longueville \cite[Chap. 3]{deLongueville2013}.
The appendices A--E to this book also provide a concise and user-friendly account of the Algebraic Topology tools
employed.
See also Miller \cite{Miller2013}.

Gorenstein \cite{Gorenstein1980} is a standard text on finite groups. 
The book by Smith \cite{Smith-SubgroupComplexes}
contains a wealth of material  on subgroup lattices and
can serve as our general reference for these.
\end{bibpar}

\begin{exs}
\item
What kind of values of~$c(\FF)$ are possible for graph properties
of graphs on $n$~vertices? For monotone properties,
it is assumed that one has $c(\FF)\in \{0,m\}$, and this is
proved if $n$ is a prime power.
In general, it is known that $c(\FF)\ge 2n-4$
unless $c(\FF)=0$, by Bollob\'as \& Eldridge \cite{BollobasEldridge},
see \cite[Sect. VIII.5]{Bollobas}.

\item
Show that the digraph property ``has a sink'' has complexity 
\[
c(\FF_{sink})\le 3(n-1)-\lfloor\log_2(n)\rfloor.
\]
Can you also prove that for any non-trivial digraph property
one has $c(\FF)\ge c(\FF_{sink})$?\\
(This is stated in Best, van Emde Boas \& Lenstra \cite[p. 17]{BEBL};
there are analogous results by Bollob\'as \& Eldridge
\cite{BollobasEldridge} \cite[Sect.~VIII.5]{Bollobas} in
a different model for digraphs.)

\item
Show that if a complex $\Delta$ corresponds to
a non-evasive monotone graph property, then
it has a complete $1$-skeleton.

\item
Give examples of simplicial complexes that are contractible, but not collapsible.
(The “dunce hat” is a key word for a search in the literature~\ldots\,)
  
\item
Assume that when testing some unknown set $A$ with respect to
a set system $\FF$, you always get the answer YES
if there is \emph{any} set $A\in\FF$ for which this YES and all the previous answers are correct,
that is, unless this “YES” would allow you to conclude $A\notin\FF$ at this point.
\begin{itemize}
\item[(i)]
Show that with this type of answers you \emph{always} need $m$ questions
for \emph{any} algorithm (and thus $\FF$ is evasive) if and only if
$\FF$ satisfies the following property:
\begin{itemize}
\item[$(*)$]
for any $e\in A\in\FF$ there is some $f\in E\sm A$ such that
$A\sm\{e\}\cup \{f\}\in\FF$.
\end{itemize}
\item[(ii)]
Show that for $n\ge 5$, the family $\FF$ of edge sets of planar graphs
satisfies property~$(*)$.
\item[(iii)]
Give other examples of graph properties that satisfy $(*)$, and 
are thus evasive.
\end{itemize}
(This is the ``simple strategy'' of Milner \& Welsh \cite{MilnerWelsh};
see Bollob\'as \cite[p. 406]{Bollobas}.)

\item 
Let $\Delta$ be a vertex-homogeneous simplicial complex with $n$ vertices
and Euler characteristic $\chi(\Delta) = -1$. Suppose that 
$n= p_1^{e_1} \cdots p_k^{e_k}$ is the prime factorization of~$n$ and
let  $m= \max \{ p_1^{e_1}, \dots , p_k^{e_k}\}$. Prove that
$\dim \Delta \ge m-1.$

\item
Let $W^q_n$ be the set of all words of length $n$ in the alphabet
$\{1, 2, \dots , q\}$, $q\ge 2$. For subsets $\FF \subseteq W^q_n$,
let $c(\FF)$ be the least number of inspections of single letters (or rather,
positions) that the best algorithm needs in the worst case $s\in W^q_n$
in order to decide whether $s\in \FF$.

Define the polynomial
\[
	p_{\FF} (x_1, \dots , x_q)=\sum_{s\in \FF} x_1^{\mu_1} \dots x_q^{\mu_q},
\]
where $\mu_i =\# \{j : s_j=i\}$ for $s=s_1 \cdots s_q$.

Show that 
\[
	(x_1 + \dots + x_q)^{n-c(\FF)}  \ \ \big|\ \     p_{\FF} (x_1, \dots , x_q).
\] 
\end{exs}
 
 
\bigskip
\section{Appendix: Fixed point theorems and homology}\label{sec5}
\subsection{Lefschetz' theorem} 

Fixed point theorems are 
``global--local tools'': From global information about a space
(such as its homology) they derive local effects, such as the 
existence of special points where ``something happens.''

Of course, in applications to combinatorial problems we need
to combine them with suitable 
``continuous--discrete tools'': From continous effects, such as
topological information about continuous maps of simplicial complexes, 
we have to find our way back to combinatorial information.

In this Appendix we assume familiarity with more Algebra and
Algebraic Topology than in other parts of these lecture notes, including
some basic finite group theory, chain complexes, etc.
As this is meant to be a reference and survey section, no detailed proofs will be given. 
A main result we head for is Oliver's theorem \ref{t:oliver},
which is needed in Section~\ref{sec:Evasiveness}.
On the way to this, skim or skip, 
depending on your tastes and familiarity\footnote{See \cite{Mat-top} for a detailed discussion of simplicial
complexes, their geometric realizations, etc. In particular,
we use the notation $\polyh K $ for the polyhedron (the geometric realization 
of a simplicial complex $\Delta$).} with these notions.
 
A powerful tool on our agenda (which yields a classical proof for
Brouwer's fixed point theorem and some of its extensions) is Hopf's trace theorem. 
For this let $V$ be any finite-dimensional vector space, 
or a free abelian group of finite rank.
When we consider an endomorphism $g\: V\to V$ 
then the \emph{trace} $\trace(g)$ is the
sum of the diagonal elements of the matrix that represents~$g$.
The trace is independent of the basis chosen for~$V$.
In the case when $V$ is a free abelian group, then $\trace(g)$ is an integer.

\begin{Theorem} [The Hopf trace theorem]\label{t:hopf-trace} 
Let $\Delta$ be a finite simplicial complex, 
let $f\: \polyh{\Delta}\to \polyh{\Delta}$ be a self-map, and denote 
by $f_{\#i}$ resp.\ $f_{*i}$ the maps that $f$ induces on 
$i$-dimensional chain groups resp.\ homology groups.

Using an arbitrary field of coefficients $\mathbf{k}$, one has
\[
\sum_i (-1)^i \trace(f_{\#i})\ \ =\ \ 
\sum_i (-1)^i \trace(f_{* i}).  
\]
The same identity holds if we use integer coefficients,
and compute the traces for homology in the quotients
$H_i(\Delta,\Z)/T_i(\Delta,\Z)$ of the homology groups modulo their torsion subgroups;
these quotients are free abelian groups.
\end{Theorem}

This theorem is remarkable as it allows to compute a topological invariant that
      depends solely on the homotopy class of~$f$, by means of a simple 
      combinatorial counting.  
The proof for this uses the definition of simplicial homology, and simple linear
algebra; we refer to Munkres \cite[Thm. 22.1]{Munkres} 
or Bredon~\cite[Sect. IV.23]{Bredon-at}.
\medskip

For an arbitrary coefficient field $\mathbf{k}$, the \emph{Lefschetz number}
of the map~$f\:\polyh{\Delta}\to \polyh{\Delta}$ is defined as
\[
L_{\mathbf{k}}(f)\assg \sum_i (-1)^i \trace(f_{* i})\ \ \in \mathbf{k}.
\]
Similarly, taking integral homology modulo torsion,  
the \emph{integral Lefschetz number} is defined as
\[
L(f)\assg\sum_i (-1)^i \trace(f_{* i})\ \ \in \Z.
\]
The universal coefficient theorems imply that one always has
$L_\Q(f)=L(f)$: Thus the integral  Lefschetz number $L(f)$ can be computed
in rational homology, but it is an integer.

The \emph{Euler characteristic} of a complex $\Delta$ coincides with the 
Lefschetz number of the identity map $\id_\Delta\: \polyh{\Delta}\to \polyh{\Delta}$,
\[
\chi(\Delta) = L(\id_\Delta),
\quad \textrm{where }
\trace((\id_\Delta)_{*i})=\beta_i(\Delta).
\]
Thus the Hopf trace theorem
yields that the Euler characteristic of a finite simplicial complex $\Delta$
can be defined resp.\ computed without a reference to homology, 
simply as the alternating sum of the face numbers
of the complex~$\Delta$, where $f_i=F_i(\Delta)$ denotes the number of $i$-dimensional
faces of~$\Delta$:
\[
\chi(\Delta)\assg f_0(\Delta) - f_1(\Delta) + f_2(\Delta) - \cdots .
\]
This is then a finite sum that ends with $(-1)^d f_d(\Delta)$ if
$\Delta$ has dimension~$d$. 
Thus the Hopf trace theorem applied to the identity 
map just reproduces the Euler--Poincar\'e formula.
This proves, for example, the $d$-dimensional Euler polyhedron formula,
not only for polytopes, but also for general spheres, shellable or not 
(as discussed in Ziegler \cite{Z-shellballs}).
The Hopf trace formula also has powerful combinatorial applications,
see Ziegler \cite{Ziegler2002}. Howver, for us its
main consequence is the following theorem, which is a vast generalization of the
Brouwer fixed point theorem.

\begin{Theorem} [The Lefschetz fixed point theorem]\label{t:LefschetzHopf}
Let $\Delta$ be a finite simplicial complex, and $\mathbf{k}$ an arbitrary field.
If a self-map $f\: \polyh{\Delta}\to \polyh{\Delta}$ 
has Lefschetz number $L_{\mathbf{k}}(f)\neq0$, then
$f$ and every map homotopic to $f$ have a fixed point.

In particular, if $\Delta$ is $\Z_p$-acyclic for 
some prime $p$, then every continuous
map $f\: \polyh{\Delta}\to \polyh{\Delta}$ has a fixed point. 
\end{Theorem}

(A complex is $\Z_p$-acyclic if its reduced homology with $\Z_p$-coefficients vanishes.
That is, in terms of homology it looks like a contractible space, say a $d$-ball.)

\proofheader{Proof \rm(Sketch)}
For a finite simplicial complex $\Delta$,
the polyhedron $\polyh{\Delta}$ is compact. So if $f$ does not have a fixed point,
there is some $\eps>0$ such that $|f(x)-x|>\eps$ for all $x\in \Delta$.
Now take a subdivision $\Delta'$ of~$\Delta$ into
simplices of diameter smaller than $\eps$, and a simplicial
approximation of error smaller than $\eps/2$, 
so that the simplicial approximation $f':\Delta'\rightarrow\Delta'$,
which is homotopic to $f$, does not have a fixed point, either.
 
Now apply the trace theorem to see that $L_{\mathbf{k}}(f)$
      is zero, contrary to the assumption, where the induced map $f'_{*0}=f_{*0}$ 
in $0$-dimensional homology is the identity.
\endproof

Note that Brouwer's fixed point theorem~\ref{t:brouwer}
is the special case of Theorem~\ref{t:LefschetzHopf}  when
$\Delta$ triangulates a ball.  

For a reasonably large class of spaces, a converse to the
Lefschetz fixed point theorem is also true: 
If $L(f)=0$, then $f$ is homotopic to a map without fixed points.
See Brown \cite[Chap. VIII]{Brown}.

\subsection{The theorems of Smith and Oliver}\label{sec4.2}

In addition to the usual game of
connections between  graphs, posets, complexes and spaces,
we will now add groups. Namely we will
discuss some useful topological effects
caused by symmetry, that is, by finite group actions. 

A (finite) group $G$ \emph{acts} on a (finite) simplicial complex
$\Delta$ if each group element corresponds to a permutation of the vertices of~$\Delta$,
where composition of group elements corresponds to composition of permutations,  
in such a way that
$g(A)\assg \{ gv\sep v\in A\}$ is a face of~$\Delta$ for all $g\in G$ and
for all $A\in \Delta$. This action on the vertices is extended to
the geometric realization of the complex~$\Delta$, so that
$G$ acts as a group of simplicial homeomorphisms
$g\: \polyh{\Delta}\to \polyh{\Delta}$.

The action is \emph{faithful} if only the identity element
in $G$ acts as the identity permutation. In general,  
the set $G_0\assg\{g\in G\sep gv=v\mbox{ for all }v\in\vert(\Delta)\}$
is a normal subgroup of~$G$. Hence we get that the quotient group
$G/G_0$ acts faithfully on~$\Delta$, and we usually 
only consider faithful actions. In this case, we can interpret
$G$ as a subgroup of the \emph{symmetry group} of the complex~$\Delta$.
The action is \emph{vertex transitive} if for any two vertices $v,w$
of~$\Delta$ there is a group element $g\in G$ with $gv=w$.

A \emph{fixed point} (also known as \emph{stable point}) 
of a group action is a point $x\in \polyh{\Delta}$ that 
satisfies $gx=x$ for all $g\in G$. We denote the
set of all fixed points by $\Delta^G$. Note that $\Delta^G$ is in
general not a subcomplex of~$\Delta$.
 
\begin{Example} 
Let $\Delta=2^{[3]}$ be the complex of a triangle, and let
$G=\Z_3$ be the cyclic group (a proper subgroup of the symmetry
group $\Symm_3$), acting such that a 
generator cyclically permutes the vertices, $1\mapsto2\mapsto3\mapsto1$. 
\[
\begin{picture}(0,0)%
\epsfig{file=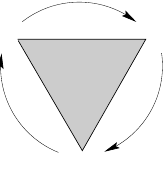}
\end{picture}%
\setlength{\unitlength}{0.00030000in}%
\begin{picture}(3616,3695)(5375,-5995)
\put(8806,-3151){\makebox(0,0)[lb]{2}}
\put(5446,-3166){\makebox(0,0)[lb]{1}}
\put(7156,-5956){\makebox(0,0)[lb]{3}}
\end{picture}
\]
This is a faithful action; its fixed point set consists of
the center of the triangle only---this is not a subcomplex of~$\Delta$,
although it corresponds to a subcomplex of the barycentric subdivision~$\sd(\Delta)$.
\end{Example}

\begin{Lemma}[Two barycentric subdivisions]\label{l:fps}~
\begin{compactenum}[\rm(1)]
 \item
After replacing $\Delta$ by its barycentric subdivision
(informally, let $\Delta :=\mathrm{sd}({\Delta})$), we 
get that the fixed point set $\Delta^G$ is a subcomplex of~$\Delta$.  
\item
After replacing $\Delta$ once again by its barycentric subdivision 
(so now $\Delta := \mathrm{sd}^2({\Delta}))$, we
even get that the quotient space $\polyh{\Delta}/G$ can be constructed from $\Delta$ by
identifying all faces with their images under the action of~$G$.
That is, the equivalence classes of faces of~$\Delta$, with the induced
partial order, form a simplicial complex that is homeomorphic
to the quotient space $\polyh{\Delta}/G$.
\end{compactenum}
\end{Lemma} 

We leave the proof as an exercise. It is not difficult; for details
and further discussion see Bredon \cite[Sect. III.1]{Bredon-tg}. 
\medskip

``Smith Theory'' was started by P.\,A. Smith \cite{Smith2} in the 
thirties.
It analyzes finite group actions on compact
spaces (such as finite simplicial complexes), providing 
relations between the structure of the group to its possible fixed point sets.
Here is one key result.

\begin{Theorem}[Smith \cite{Smith3}]\label{t:smith} 
If $P$ is a $p$-group (that is, a finite group of order $|P|=p^t$ for a prime $p$ and some $t>0$),
acting on a complex $\Delta$ that is $\Z_p$-acyclic,
then the fixed point set $\Delta^P$ is $\Z_p$-acyclic as well.
In particular, it is not empty.
\end{Theorem}

\proofheader{Proof \rm(Sketch)}
The key is that, with the preparations of Lemma~\ref{l:fps}, the maps that
$f$ induces on the chain groups (with $\Z_p$ coefficients)
nicely restrict to the chain groups on the fixed point set $\Delta^P$.
Passing to traces and using the Hopf trace theorem, one can
derive that $\Delta^P$ is non-empty.
A more detailed analysis leads to the
``transfer isomorphism'' in homology, which 
proves that $\Delta^P$ must be acyclic. 

See Bredon \cite[Thm. III.5.2]{Bredon-tg} and Oliver \cite[p. 157]{Oliver},
and also de Longueville \cite[Appendix D and~E]{deLongueville2013}.
\endproof
    
On the combinatorial side, one has an Euler characteristic relation 
due to Floyd \cite{Floyd} \cite[Sect. III.4]{Bredon-tg}:
\[
\chi(\Delta)\ +\ (p-1)\chi(\Delta^{\Z_p})\ \ =\ \ p\,\chi(\Delta/{\Z_p}).
\]
If $P$ is a $p$-group (in particular for $P=\Z_p$), then this implies that
\[
\chi(\Delta^P)\equiv \chi(\Delta) \pmod p,
\]
using induction on~$t$, where $|P|=p^t$.

\begin{Theorem}[{Oliver \cite[Lemma I]{Oliver}}]\label{t:cyclic-oliver}
If $G=\Z_n$ is a cyclic group, acting on a 
$\Q$-acyclic complex $\Delta$, then the action has a fixed point.

In this case the fixed point set $\Delta^G$ has 
the Euler characteristic of a point, ${\chi}(\Delta^G)=1$.
\end{Theorem}

\proof
The first statement follows directly from the Lefschetz
fixed point theorem: Any cyclic group is generated by a single
element $g$, this element has a fixed point, this fixed point
of~$g$ is also a fixed point of all powers of~$g$, and hence of the whole
group $G$.

For the second part, take $p^t$ to be a maximal prime power that divides $n$,
consider the corresponding subgroup isomorphic to $\Z_{p^t}$,
and use induction on~$t$ and the transfer homomorphism, 
as for the previous proof.
\endproof 
 
Unfortunately, results like these may give an overly optimistic impression of
the generality of
fixed point theorems for acyclic complexes. There are fixed point free 
finite group actions on balls: Examples were constructed by
Floyd \& Richardson and others; see Bredon \cite[Sect. I.8]{Bredon-tg}.

On the positive side we have the following result due to Oliver,
which plays a central role in Section \ref{sec4}.
 
\begin{Theorem}[Oliver's Theorem I~{\cite[Prop. I]{Oliver}}]\label{t:oliver}
If $G$ has a normal subgroup $P\triangleleft G$
that is a $p$-group, such that the quotient $G/P$ is cyclic, 
acting on a complex $\Delta$ that is $\Z_p$-acyclic,
then the fixed point set $\Delta^G$ is $\Z_p$-acyclic as well.
In particular, it is not empty.
\end{Theorem}

This is as much as we will need in this chapter.
Oliver proved, in fact, a more general and complete
theorem that includes a converse.

\begin{Theorem}[Oliver's Theorem II~{\cite{Oliver}}]\label{t:oliver2}
Let $G$ be a finite group. Every action of~$G$ on a 
$\Z_p$-acyclic complex $\Delta$ has a fixed point if and only if 
$G$ has the following structure:
\begin{quote}
$G$ has normal subgroups $P\triangleleft Q\triangleleft G$
such that $P$ is a $p$-group, $G/Q$ is a $q$-group 
(for a prime $q$ that need not be distinct from $p$),
and the quotient $Q/P$ is cyclic. 
\end{quote}
In this situation one always has $\chi(\Delta^G)\equiv 1 \bmod q$.
\end{Theorem}
 
\begin{bibpar}
The Lefschetz--Hopf fixed point theorem was announced by
Lefschetz for a restriced class of complexes in 1923, with details
appearing three years later. The first proof for the general version
was by Hopf in 1929. There are generalizations, for example to
Absolute Neighborhood Retracts; see Bredon \cite[Cor. IV.23.5]{Bredon-at} 
and Brown \cite[Chap. IIII]{Brown}. We refer to Brown's book~\cite{Brown}.

We refer to Bredon \cite[Chapter III]{Bredon-tg} for a nice 
textbook treatment of Smith Theory.
The book by de Longueville \cite[Appendix E]{deLongueville2013} 
also has a very accessible discussion
of the fixed point theorems of Smith and Oliver. The exercises 
concerning fixed point sets of poset maps $P\rightarrow P$ are drawn from 
Baclawski \& Bj\"orner \cite{BB79}.
\end{bibpar}

\begin{exs}
\item
Verify directly that if $f$ maps $\polyh{T}$ to $\polyh{T}$, where $T$ is a 
      graph-theoretic tree, then $f$ has a fixed point. 
	  
How would you derive this from the Lefschetz fixed point theorem?
	  
\item Let $P$ be a poset (finite partially ordered set), and denote by
$\Delta(P)$ its order complex (whose faces are the totally ordered 
subsets).  
Suppose that $f\: P \rightarrow P$ is an order-preserving mapping
with fixed point set $P^f := \{x\in P \mid f(x)=x\}$. \\
(a) Show that if
$\Delta(P)$ is acyclic over some field, then
\[
	\mu(P^f)=0,
\]
where $\mu(P^f)$ denotes the \emph{M\"obius function} (reduced 
Euler characteristic) of $\Delta(P^f)$. In particular, $P^f$
is not empty. \\
(b) Does it follow also that $P^f$ itself is acyclic?
\item 
Suppose now that $f\: P \rightarrow P$ is order-reversing 
and let  $P_f := \{x\in P \mid x=f^2(x)\le f(x)\}$. Show that if
$\Delta(P)$ is acyclic over some field, then
\[
	\mu(P_f)=0.
\] 
In particular, if $f$ has no fixed edge  
(i.e., no $x$ such  that $x=f^2(x)< f(x)$)
then $f$ has a unique fixed point.
\end{exs}

\subsubsection*{Acknowledgements}
We are grateful to Marie-Sophie Litz and to the referees for very careful reading
and a great number of very valuable comments and suggestions on the manuscript.
Thanks to Moritz Firsching and Stephen D. Smith, 
and in particular to Penny Haxell,
for additional references and very helpful explanations.

\begin{small}

\end{small} 
\end{document}

%% file: EPS/brouwer2-3a.pstex_t
\begin{picture}(0,0)%
\includegraphics{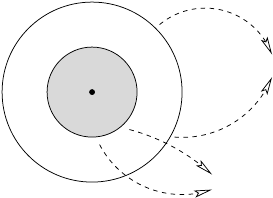}
\end{picture}%
\setlength{\unitlength}{1579sp}%
\begingroup\makeatletter\ifx\SetFigFont\undefined%
\gdef\SetFigFont#1#2#3#4#5{%
  \reset@font\fontsize{#1}{#2pt}%
  \fontfamily{#3}\fontseries{#4}\fontshape{#5}%
  \selectfont}%
\fi\endgroup%
\begin{picture}(5430,3922)(2986,-6069)
\put(7351,-5911){\makebox(0,0)[lb]{\smash{\SetFigFont{10}{12.0}{\familydefault}{\mddefault}{\updefault}$\xx\longmapsto \xx_0$}}}
\put(7951,-3511){\makebox(0,0)[lb]{\smash{\SetFigFont{10}{12.0}{\familydefault}{\mddefault}{\updefault}$\xx\longmapsto f(\xx)$}}}
\end{picture}

%% file: EPS/hexstrategy3corr.pstex_t
\begin{picture}(0,0)%
\includegraphics{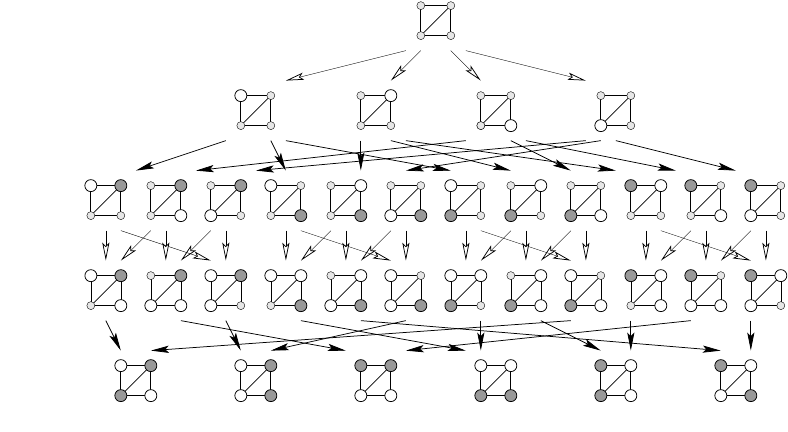}
\end{picture}%
\setlength{\unitlength}{1579sp}%
\begingroup\makeatletter\ifx\SetFigFont\undefined%
\gdef\SetFigFont#1#2#3#4#5{%
  \reset@font\fontsize{#1}{#2pt}%
  \fontfamily{#3}\fontseries{#4}\fontshape{#5}%
  \selectfont}%
\fi\endgroup%
\begin{picture}(15750,8601)(-614,-8875)
\put(1201,-1411){\makebox(0,0)[lb]{\smash{{\SetFigFont{10}{12.0}{\familydefault}{\mddefault}{\updefault}$W$ moves:}}}}
\put(1201,-3286){\makebox(0,0)[lb]{\smash{{\SetFigFont{10}{12.0}{\familydefault}{\mddefault}{\updefault}$B$ moves:}}}}
\put(-374,-7036){\makebox(0,0)[lb]{\smash{{\SetFigFont{10}{12.0}{\familydefault}{\mddefault}{\updefault}$B$ moves:}}}}
\put(-599,-5236){\makebox(0,0)[lb]{\smash{{\SetFigFont{10}{12.0}{\familydefault}{\mddefault}{\updefault}$W$ moves:}}}}
\put(1651,-8761){\makebox(0,0)[lb]{\smash{{\SetFigFont{10}{12.0}{\familydefault}{\mddefault}{\updefault}$B$ wins.}}}}
\put(4051,-8761){\makebox(0,0)[lb]{\smash{{\SetFigFont{10}{12.0}{\familydefault}{\mddefault}{\updefault}$B$ wins.}}}}
\put(11251,-8761){\makebox(0,0)[lb]{\smash{{\SetFigFont{10}{12.0}{\familydefault}{\mddefault}{\updefault}$B$ wins.}}}}
\put(8851,-8761){\makebox(0,0)[lb]{\smash{{\SetFigFont{10}{12.0}{\familydefault}{\mddefault}{\updefault}$W$ wins.}}}}
\put(13651,-8761){\makebox(0,0)[lb]{\smash{{\SetFigFont{10}{12.0}{\familydefault}{\mddefault}{\updefault}$W$ wins.}}}}
\put(6451,-8761){\makebox(0,0)[lb]{\smash{{\SetFigFont{10}{12.0}{\familydefault}{\mddefault}{\updefault}$W$ wins.}}}}
\end{picture}%

%% file: EPS/scorpion1.pstex_t
\begin{picture}(0,0)%
\includegraphics{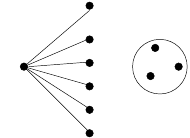}
\end{picture}%
\setlength{\unitlength}{1973sp}%
\begingroup\makeatletter\ifx\SetFigFont\undefined%
\gdef\SetFigFont#1#2#3#4#5{%
  \reset@font\fontsize{#1}{#2pt}%
  \fontfamily{#3}\fontseries{#4}\fontshape{#5}%
  \selectfont}%
\fi\endgroup%
\begin{picture}(2995,2174)(5176,-3893)
\put(5176,-3136){\makebox(0,0)[lb]{\smash{\SetFigFont{10}{12.0}{\familydefault}{\mddefault}{\updefault}$k$}}}
\end{picture}

%% file: EPS/scorpion.pstex_t
\begin{picture}(0,0)%
\includegraphics{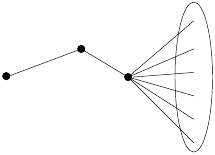}
\end{picture}%
\setlength{\unitlength}{1973sp}%
\begingroup\makeatletter\ifx\SetFigFont\undefined%
\gdef\SetFigFont#1#2#3#4#5{%
  \reset@font\fontsize{#1}{#2pt}%
  \fontfamily{#3}\fontseries{#4}\fontshape{#5}%
  \selectfont}%
\fi\endgroup%
\begin{picture}(3383,2414)(3526,-3968)
\put(3526,-2611){\makebox(0,0)[lb]{\smash{\SetFigFont{10}{12.0}{\familydefault}{\mddefault}{\updefault}1}}}
\put(4726,-2161){\makebox(0,0)[lb]{\smash{\SetFigFont{10}{12.0}{\familydefault}{\mddefault}{\updefault}2}}}
\put(5026,-3136){\makebox(0,0)[lb]{\smash{\SetFigFont{10}{12.0}{\familydefault}{\mddefault}{\updefault}$n-2$}}}
\end{picture}

%% file: EPS/inner2.pstex_t
\begin{picture}(0,0)%
\includegraphics{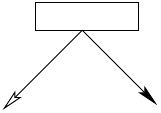}
\end{picture}%
\setlength{\unitlength}{1973sp}%
\begingroup\makeatletter\ifx\SetFigFont\undefined%
\gdef\SetFigFont#1#2#3#4#5{%
  \reset@font\fontsize{#1}{#2pt}%
  \fontfamily{#3}\fontseries{#4}\fontshape{#5}%
  \selectfont}%
\fi\endgroup%
\begin{picture}(2519,1759)(3054,-4508)
\put(3901,-3061){\makebox(0,0)[lb]{\smash{\SetFigFont{10}{12.0}{\familydefault}{\mddefault}{\updefault}$e\in A?$}}}
\put(5101,-3811){\makebox(0,0)[lb]{\smash{\SetFigFont{10}{12.0}{\familydefault}{\mddefault}{\updefault}YES}}}
\put(3226,-3811){\makebox(0,0)[lb]{\smash{\SetFigFont{10}{12.0}{\familydefault}{\mddefault}{\updefault}NO}}}
\end{picture}

%% file: EPS/k3d.pstex_t
\begin{picture}(0,0)%
\includegraphics{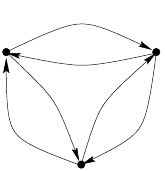}
\end{picture}%
\setlength{\unitlength}{1973sp}%
\begin{picture}(2536,2677)(4733,-4628)
\put(6421,-3811){\makebox(0,0)[lb]{{32}}}
\put(7021,-4096){\makebox(0,0)[lb]{{23}}}
\put(4936,-3886){\makebox(0,0)[lb]{{31}}}
\put(5611,-3496){\makebox(0,0)[lb]{{13}}}
\put(5881,-2821){\makebox(0,0)[lb]{{21}}}
\put(5866,-2161){\makebox(0,0)[lb]{{12}}}
\end{picture}

%% file: EPS/tree2.pstex_t
\begin{picture}(0,0)%
\includegraphics{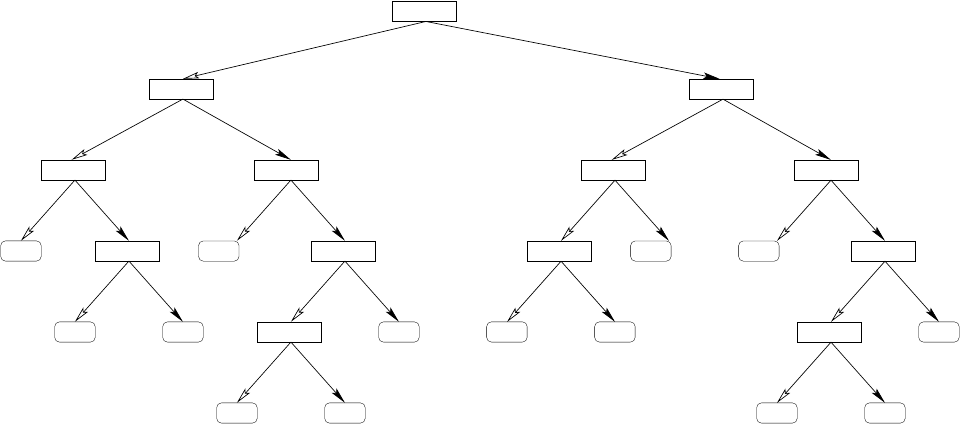}
\end{picture}%
\setlength{\unitlength}{1421sp}%
\begingroup\makeatletter\ifx\SetFigFont\undefined%
\gdef\SetFigFont#1#2#3#4#5{%
  \reset@font\fontsize{#1}{#2pt}%
  \fontfamily{#3}\fontseries{#4}\fontshape{#5}%
  \selectfont}%
\fi\endgroup%
\begin{picture}(21324,9398)(2690,-9224)
\put(13576,-7262){\makebox(0,0)[lb]{\smash{\SetFigFont{7}{8.4}{\familydefault}{\mddefault}{\updefault}NO}}}
\put(15901,-7262){\makebox(0,0)[lb]{\smash{\SetFigFont{7}{8.4}{\familydefault}{\mddefault}{\updefault}YES}}}
\put(11177,-7262){\makebox(0,0)[lb]{\smash{\SetFigFont{7}{8.4}{\familydefault}{\mddefault}{\updefault}NO}}}
\put(9977,-9062){\makebox(0,0)[lb]{\smash{\SetFigFont{7}{8.4}{\familydefault}{\mddefault}{\updefault}NO}}}
\put(21977,-9062){\makebox(0,0)[lb]{\smash{\SetFigFont{7}{8.4}{\familydefault}{\mddefault}{\updefault}NO}}}
\put(23177,-7262){\makebox(0,0)[lb]{\smash{\SetFigFont{7}{8.4}{\familydefault}{\mddefault}{\updefault}NO}}}
\put(19502,-9062){\makebox(0,0)[lb]{\smash{\SetFigFont{7}{8.4}{\familydefault}{\mddefault}{\updefault}YES}}}
\put(7502,-9062){\makebox(0,0)[lb]{\smash{\SetFigFont{7}{8.4}{\familydefault}{\mddefault}{\updefault}YES}}}
\put(6302,-7262){\makebox(0,0)[lb]{\smash{\SetFigFont{7}{8.4}{\familydefault}{\mddefault}{\updefault}YES}}}
\put(3977,-7262){\makebox(0,0)[lb]{\smash{\SetFigFont{7}{8.4}{\familydefault}{\mddefault}{\updefault}NO}}}
\put(2777,-5462){\makebox(0,0)[lb]{\smash{\SetFigFont{7}{8.4}{\familydefault}{\mddefault}{\updefault}NO}}}
\put(4953,-5462){\makebox(0,0)[lb]{\smash{\SetFigFont{7}{8.4}{\familydefault}{\mddefault}{\updefault}$31\in A?$}}}
\put(3826,-3286){\makebox(0,0)[lb]{\smash{\SetFigFont{7}{8.4}{\familydefault}{\mddefault}{\updefault}1}}}
\put(8551,-3286){\makebox(0,0)[lb]{\smash{\SetFigFont{7}{8.4}{\familydefault}{\mddefault}{\updefault}3}}}
\put(7353,-5462){\makebox(0,0)[lb]{\smash{\SetFigFont{7}{8.4}{\familydefault}{\mddefault}{\updefault}NO}}}
\put(9753,-5462){\makebox(0,0)[lb]{\smash{\SetFigFont{7}{8.4}{\familydefault}{\mddefault}{\updefault}$31\in A?$}}}
\put(8553,-7262){\makebox(0,0)[lb]{\smash{\SetFigFont{7}{8.4}{\familydefault}{\mddefault}{\updefault}$32\in A?$}}}
\put(15755,-3665){\makebox(0,0)[lb]{\smash{\SetFigFont{7}{8.4}{\familydefault}{\mddefault}{\updefault}$21\in A?$}}}
\put(20553,-3286){\makebox(0,0)[lb]{\smash{\SetFigFont{7}{8.4}{\familydefault}{\mddefault}{\updefault}3}}}
\put(19355,-5465){\makebox(0,0)[lb]{\smash{\SetFigFont{7}{8.4}{\familydefault}{\mddefault}{\updefault}NO}}}
\put(21755,-5465){\makebox(0,0)[lb]{\smash{\SetFigFont{7}{8.4}{\familydefault}{\mddefault}{\updefault}$31\in A?$}}}
\put(20555,-7262){\makebox(0,0)[lb]{\smash{\SetFigFont{7}{8.4}{\familydefault}{\mddefault}{\updefault}$32\in A?$}}}
\put(18155,-1865){\makebox(0,0)[lb]{\smash{\SetFigFont{7}{8.4}{\familydefault}{\mddefault}{\updefault}$23\in A?$}}}
\put(20480,-3665){\makebox(0,0)[lb]{\smash{\SetFigFont{7}{8.4}{\familydefault}{\mddefault}{\updefault}$13\in A?$}}}
\put(15828,-3288){\makebox(0,0)[lb]{\smash{\SetFigFont{7}{8.4}{\familydefault}{\mddefault}{\updefault}2}}}
\put(14556,-5466){\makebox(0,0)[lb]{\smash{\SetFigFont{7}{8.4}{\familydefault}{\mddefault}{\updefault}$32\in A?$}}}
\put(16956,-5466){\makebox(0,0)[lb]{\smash{\SetFigFont{7}{8.4}{\familydefault}{\mddefault}{\updefault}NO}}}
\put(3753,-3663){\makebox(0,0)[lb]{\smash{\SetFigFont{7}{8.4}{\familydefault}{\mddefault}{\updefault}$21\in A?$}}}
\put(6153,-1865){\makebox(0,0)[lb]{\smash{\SetFigFont{7}{8.4}{\familydefault}{\mddefault}{\updefault}$13\in A?$}}}
\put(8478,-3665){\makebox(0,0)[lb]{\smash{\SetFigFont{7}{8.4}{\familydefault}{\mddefault}{\updefault}$23\in A?$}}}
\put(11553,-138){\makebox(0,0)[lb]{\smash{\SetFigFont{7}{8.4}{\familydefault}{\mddefault}{\updefault}$12\in A?$}}}
\put(9526,-1186){\makebox(0,0)[lb]{\smash{\SetFigFont{7}{8.4}{\familydefault}{\mddefault}{\updefault}NO}}}
\put(3826,-4636){\makebox(0,0)[lb]{\smash{\SetFigFont{7}{8.4}{\familydefault}{\mddefault}{\updefault}NO}}}
\put(10000,-4636){\makebox(0,0)[lb]{\smash{\SetFigFont{7}{8.4}{\familydefault}{\mddefault}{\updefault}YES}}}
\put(8626,-8236){\makebox(0,0)[lb]{\smash{\SetFigFont{7}{8.4}{\familydefault}{\mddefault}{\updefault}NO}}}
\put(9826,-6436){\makebox(0,0)[lb]{\smash{\SetFigFont{7}{8.4}{\familydefault}{\mddefault}{\updefault}NO}}}
\put(8626,-4636){\makebox(0,0)[lb]{\smash{\SetFigFont{7}{8.4}{\familydefault}{\mddefault}{\updefault}NO}}}
\put(15826,-4636){\makebox(0,0)[lb]{\smash{\SetFigFont{7}{8.4}{\familydefault}{\mddefault}{\updefault}NO}}}
\put(14626,-6436){\makebox(0,0)[lb]{\smash{\SetFigFont{7}{8.4}{\familydefault}{\mddefault}{\updefault}NO}}}
\put(6301,-6436){\makebox(0,0)[lb]{\smash{\SetFigFont{7}{8.4}{\familydefault}{\mddefault}{\updefault}YES}}}
\put(11101,-6436){\makebox(0,0)[lb]{\smash{\SetFigFont{7}{8.4}{\familydefault}{\mddefault}{\updefault}YES}}}
\put(15901,-6436){\makebox(0,0)[lb]{\smash{\SetFigFont{7}{8.4}{\familydefault}{\mddefault}{\updefault}YES}}}
\put(17101,-4636){\makebox(0,0)[lb]{\smash{\SetFigFont{7}{8.4}{\familydefault}{\mddefault}{\updefault}YES}}}
\put(16276,-1036){\makebox(0,0)[lb]{\smash{\SetFigFont{7}{8.4}{\familydefault}{\mddefault}{\updefault}YES}}}
\put(8251,-2836){\makebox(0,0)[lb]{\smash{\SetFigFont{7}{8.4}{\familydefault}{\mddefault}{\updefault}YES}}}
\put(5101,-4636){\makebox(0,0)[lb]{\smash{\SetFigFont{7}{8.4}{\familydefault}{\mddefault}{\updefault}YES}}}
\put(5026,-6436){\makebox(0,0)[lb]{\smash{\SetFigFont{7}{8.4}{\familydefault}{\mddefault}{\updefault}NO}}}
\put(5626,-2836){\makebox(0,0)[lb]{\smash{\SetFigFont{7}{8.4}{\familydefault}{\mddefault}{\updefault}NO}}}
\put(20176,-2836){\makebox(0,0)[lb]{\smash{\SetFigFont{7}{8.4}{\familydefault}{\mddefault}{\updefault}YES}}}
\put(21901,-8236){\makebox(0,0)[lb]{\smash{\SetFigFont{7}{8.4}{\familydefault}{\mddefault}{\updefault}YES}}}
\put(20626,-4636){\makebox(0,0)[lb]{\smash{\SetFigFont{7}{8.4}{\familydefault}{\mddefault}{\updefault}NO}}}
\put(21826,-6436){\makebox(0,0)[lb]{\smash{\SetFigFont{7}{8.4}{\familydefault}{\mddefault}{\updefault}NO}}}
\put(21901,-4636){\makebox(0,0)[lb]{\smash{\SetFigFont{7}{8.4}{\familydefault}{\mddefault}{\updefault}YES}}}
\put(23101,-6436){\makebox(0,0)[lb]{\smash{\SetFigFont{7}{8.4}{\familydefault}{\mddefault}{\updefault}YES}}}
\put(20626,-8236){\makebox(0,0)[lb]{\smash{\SetFigFont{7}{8.4}{\familydefault}{\mddefault}{\updefault}NO}}}
\put(9901,-8236){\makebox(0,0)[lb]{\smash{\SetFigFont{7}{8.4}{\familydefault}{\mddefault}{\updefault}YES}}}
\put(17701,-2836){\makebox(0,0)[lb]{\smash{\SetFigFont{7}{8.4}{\familydefault}{\mddefault}{\updefault}NO}}}
\end{picture}

%% file: EPS/boolean5.pstex_t
\begin{picture}(0,0)%
\includegraphics{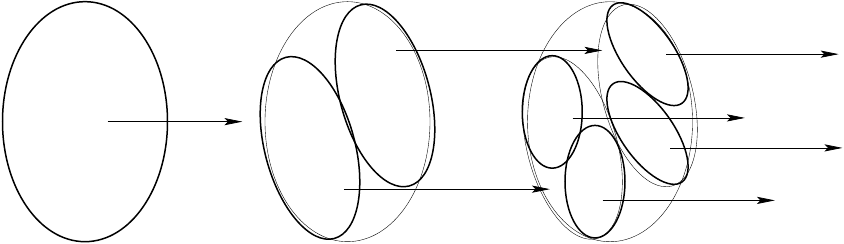}
\end{picture}%
\setlength{\unitlength}{1579sp}%
\begingroup\makeatletter\ifx\SetFigFont\undefined%
\gdef\SetFigFont#1#2#3#4#5{%
  \reset@font\fontsize{#1}{#2pt}%
  \fontfamily{#3}\fontseries{#4}\fontshape{#5}%
  \selectfont}%
\fi\endgroup%
\begin{picture}(16845,4844)(5528,-5183)
\put(14551,-1186){\makebox(0,0)[lb]{\smash{\SetFigFont{8}{9.6}{\familydefault}{\mddefault}{\updefault}$f\in A?$}}}
\put(14476,-3961){\makebox(0,0)[lb]{\smash{\SetFigFont{8}{9.6}{\familydefault}{\mddefault}{\updefault}$g\in A?$}}}
\put(9076,-2611){\makebox(0,0)[lb]{\smash{\SetFigFont{8}{9.6}{\familydefault}{\mddefault}{\updefault}$e\in A?$}}}
\end{picture}

%% file: EPS/orbit.pstex_t
\begin{picture}(0,0)%
\includegraphics{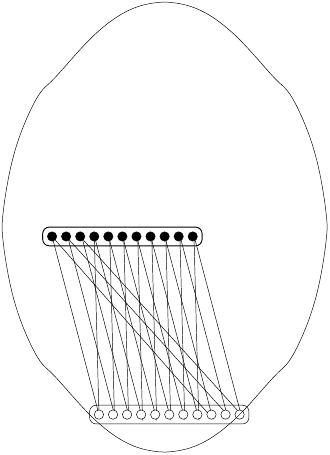}
\end{picture}%
\setlength{\unitlength}{1973sp}%
\begingroup\makeatletter\ifx\SetFigFont\undefined%
\gdef\SetFigFont#1#2#3#4#5{%
  \reset@font\fontsize{#1}{#2pt}%
  \fontfamily{#3}\fontseries{#4}\fontshape{#5}%
  \selectfont}%
\fi\endgroup%
\begin{picture}(5224,7224)(12989,-7573)
\put(13276,-1561){\makebox(0,0)[lb]{\smash{\SetFigFont{10}{12.0}{\familydefault}{\mddefault}{\updefault}$2^E$}}}
\put(13801,-3811){\makebox(0,0)[lb]{\smash{\SetFigFont{10}{12.0}{\familydefault}{\mddefault}{\updefault}$\OO$}}}
\put(17101,-7111){\makebox(0,0)[lb]{\smash{\SetFigFont{10}{12.0}{\familydefault}{\mddefault}{\updefault}$E$: has $p^t$ elements}}}
\end{picture}

%% file: EPS/illies3.pstex_t
\begin{picture}(0,0)%
\includegraphics{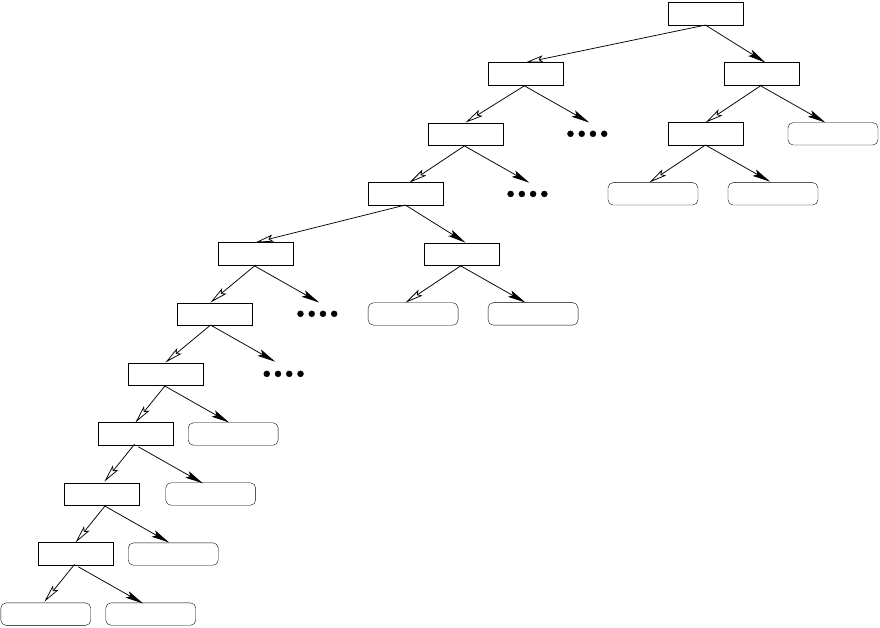}
\end{picture}%
\setlength{\unitlength}{1579sp}%
\begingroup\makeatletter\ifx\SetFigFont\undefined%
\gdef\SetFigFont#1#2#3#4#5{%
  \reset@font\fontsize{#1}{#2pt}%
  \fontfamily{#3}\fontseries{#4}\fontshape{#5}%
  \selectfont}%
\fi\endgroup%
\begin{picture}(17567,12480)(297,-12305)
\put(10204,-6139){\makebox(0,0)[lb]{\smash{\SetFigFont{8}{9.6}{\familydefault}{\mddefault}{\updefault}YES(10)}}}
\put(4804,-4939){\makebox(0,0)[lb]{\smash{\SetFigFont{8}{9.6}{\familydefault}{\mddefault}{\updefault}$11\in A?$}}}
\put(2401,-8536){\makebox(0,0)[lb]{\smash{\SetFigFont{8}{9.6}{\familydefault}{\mddefault}{\updefault}$10\in A?$}}}
\put(1732,-9742){\makebox(0,0)[lb]{\smash{\SetFigFont{8}{9.6}{\familydefault}{\mddefault}{\updefault}$6\in A?$}}}
\put(4207,-8542){\makebox(0,0)[lb]{\smash{\SetFigFont{8}{9.6}{\familydefault}{\mddefault}{\updefault}YES(10)}}}
\put(3757,-9742){\makebox(0,0)[lb]{\smash{\SetFigFont{8}{9.6}{\familydefault}{\mddefault}{\updefault}YES(6)}}}
\put(3008,-10943){\makebox(0,0)[lb]{\smash{\SetFigFont{8}{9.6}{\familydefault}{\mddefault}{\updefault}YES(9)}}}
\put(1202,-10937){\makebox(0,0)[lb]{\smash{\SetFigFont{8}{9.6}{\familydefault}{\mddefault}{\updefault}$9\in A?$}}}
\put(459,-12144){\makebox(0,0)[lb]{\smash{\SetFigFont{8}{9.6}{\familydefault}{\mddefault}{\updefault}YES(5,8)}}}
\put(2558,-12143){\makebox(0,0)[lb]{\smash{\SetFigFont{8}{9.6}{\familydefault}{\mddefault}{\updefault}YES(5)}}}
\put(3983,-6143){\makebox(0,0)[lb]{\smash{\SetFigFont{8}{9.6}{\familydefault}{\mddefault}{\updefault}$12\in A?$}}}
\put(3005,-7340){\makebox(0,0)[lb]{\smash{\SetFigFont{8}{9.6}{\familydefault}{\mddefault}{\updefault}$7\in A?$}}}
\put(7806,-6141){\makebox(0,0)[lb]{\smash{\SetFigFont{8}{9.6}{\familydefault}{\mddefault}{\updefault}YES(8,12)}}}
\put(7804,-3769){\makebox(0,0)[lb]{\smash{\SetFigFont{8}{9.6}{\familydefault}{\mddefault}{\updefault}$4\in A?$}}}
\put(8931,-4971){\makebox(0,0)[lb]{\smash{\SetFigFont{8}{9.6}{\familydefault}{\mddefault}{\updefault}$7\in A?$}}}
\put(9007,-2572){\makebox(0,0)[lb]{\smash{\SetFigFont{8}{9.6}{\familydefault}{\mddefault}{\updefault}$3\in A?$}}}
\put(12603,-3753){\makebox(0,0)[lb]{\smash{\SetFigFont{8}{9.6}{\familydefault}{\mddefault}{\updefault}YES(5,9)}}}
\put(13789,-2569){\makebox(0,0)[lb]{\smash{\SetFigFont{8}{9.6}{\familydefault}{\mddefault}{\updefault}$10\in A?$}}}
\put(16127,-2567){\makebox(0,0)[lb]{\smash{\SetFigFont{8}{9.6}{\familydefault}{\mddefault}{\updefault}YES(7,10)}}}
\put(15018,-3753){\makebox(0,0)[lb]{\smash{\SetFigFont{8}{9.6}{\familydefault}{\mddefault}{\updefault}YES(7)}}}
\put(13787,-182){\makebox(0,0)[lb]{\smash{\SetFigFont{8}{9.6}{\familydefault}{\mddefault}{\updefault}$1\in A?$}}}
\put(10174,-1369){\makebox(0,0)[lb]{\smash{\SetFigFont{8}{9.6}{\familydefault}{\mddefault}{\updefault}$2\in A?$}}}
\put(14914,-1369){\makebox(0,0)[lb]{\smash{\SetFigFont{8}{9.6}{\familydefault}{\mddefault}{\updefault}$4\in A?$}}}
\end{picture}

%% file: EPS/Kkn.pstex_t
\begin{picture}(0,0)%
\includegraphics{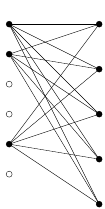}
\end{picture}%
\setlength{\unitlength}{1579sp}%
\begingroup\makeatletter\ifx\SetFigFont\undefined%
\gdef\SetFigFont#1#2#3#4#5{%
  \reset@font\fontsize{#1}{#2pt}%
  \fontfamily{#3}\fontseries{#4}\fontshape{#5}%
  \selectfont}%
\fi\endgroup%
\begin{picture}(2018,4132)(3451,-5228)
\put(3451,-1336){\makebox(0,0)[lb]{\smash{\SetFigFont{8}{9.6}{\familydefault}{\mddefault}{\updefault}$V$}}}
\put(5251,-1336){\makebox(0,0)[lb]{\smash{\SetFigFont{8}{9.6}{\familydefault}{\mddefault}{\updefault}$W$}}}
\end{picture}